\documentclass[draft,10pt]{article}
\usepackage{amssymb}
\usepackage{verbatim,cite}
\usepackage{latexsym}
\usepackage{amsmath}
\usepackage{eucal}
\newtheorem{theorem}{Theorem}[section]

\newtheorem{sub}{\name}[section]
\newcommand{\bs}{
\begin{sub}}
\newcommand{\es}{
\end{sub}}
\newcommand{\bsl}[1]{
\begin{sub}\label{#1}}
\newcommand{\bth}[1]{\def\name{Theorem}
\begin{sub}\label{t:#1}}
\newcommand{\blemma}[1]{\def\name{Lemma}
\begin{sub}\label{l:#1}}
\newcommand{\bcor}[1]{\def\name{Corollary}
\begin{sub}\label{c:#1}}
\newcommand{\bdef}[1]{\def\name{Definition}
\begin{sub}\label{d:#1}}
\newcommand{\bprop}[1]{\def\name{Proposition}
\begin{sub}\label{p:#1}}
\newcommand{\brem}[1]{\def\name{Remark}
\begin{sub}\label{r:#1}}

\newcommand{\rth}[1]{Theorem~\ref{t:#1}}

\newcommand{\rprop}[1]{Proposition~\ref{p:#1}}


\newcommand{\BA}{
\begin{array}}
\newcommand{\EA}{
\end{array}}
\newcommand{\la}{\label}
\newcommand{\BAN}{\renewcommand{\arraystretch}{1.2}
\setlength{\arraycolsep}{2pt}
\begin{array}}
\newcommand{\BAV}[2]{\renewcommand{\arraystretch}{#1}
\setlength{\arraycolsep}{#2}
\begin{array}}
\newcommand{\BSA}{
\begin{subarray}}
\newcommand{\ESA}{
\end{subarray}}
\newcommand{\BAL}{
\begin{aligned}}
\newcommand{\EAL}{
\end{aligned}}
\newcommand{\BALG}{
\begin{alignat}}
\newcommand{\EALG}{
\end{alignat}}
\newcommand{\BALGN}{
\begin{alignat*}}
\newcommand{\EALGN}{
\end{alignat*}}
\newcommand{\note}[1]{\textit{#1.}\hspace{2mm}}
\newcommand{\Proof}{\note{Proof}}
\newcommand{\qeda}{\hspace{10mm}\hfill $\square$}

\newcommand{\Remark}{\note{Remark}}

\newcommand{\forevery}{\quad \forall}

\newcommand{\warrow}{\rightharpoonup}
\newcommand{
\paran}[1]{\left (#1 \right )}
\newcommand{\sqbr}[1]{\left [#1 \right ]}
\newcommand{\curlybr}[1]{\left \{#1 \right \}}
\newcommand{\abs}[1]{\left |#1\right |}

\def\angb<#1>{\langle #1 \rangle}

\newcommand{\opname}[1]{\mbox{\rm #1}\,}
\newcommand{\supp}{\opname{supp}}
\newcommand{\dist}{\opname{dist}}
\newcommand{\myfrac}[2]{{\displaystyle \frac{#1}{#2} }}
\newcommand{\myint}[2]{{\displaystyle \int_{#1}^{#2}}}
\newcommand {\dint}{{\displaystyle \int\!\!\int}}
\newcommand{\q}{\quad}

\newcommand{\ity}{\infty}
\newcommand{\prt}{
\partial}
\newcommand{\sms}{\setminus}

\newcommand{\ti}{\times}

\newcommand{\sbs}{\subset}

\newcommand{\ovl}{\overline}

\def\ga{\alpha}     \def\gb{\beta}       \def\gga{\gamma}
       \def\gd{\delta}      \def\ge{\epsilon}
\def\gth{\theta}                         \def\vge{\varepsilon}
\def\gf{\phi}           
      \def\gk{\kappa}      
\def\gm{\mu}                 \def\gp{\pi}
    \def\gr{\rho}        
\def\gs{\sigma}       \def\gt{\tau}
      \def\gw{\omega}
                
     \def\Gd{\Delta}

\def\Gw{\Omega}              

   \def\BBN {\mathbb N}    
   \def\BBR {\mathbb R}

\allowdisplaybreaks
\begin{document}
\title {\bf  The balance between diffusion and absorption in semilinear parabolic equations\footnote {To appear in{ \it Rend. Lincei, Mat. Appl.}}}
\author{{\bf\large Andrey Shishkov}
\hspace{2mm}\vspace{3mm}\\
{\it \normalsize  Institute of Applied Mathematics and Mechanics of NAS of Ukraine},\\
{\it\normalsize R. Luxemburg str. 74, 83114 Donetsk, Ukraine}\\
\vspace{3mm}\\
{\bf\large Laurent V\'eron}
\vspace{3mm}\\
{\it\normalsize Laboratoire de Math\'ematiques et Physique Th\'eorique, CNRS UMR 6083},
\\
{\it\normalsize Universit\'e Fran\c{c}ois-Rabelais,  37200 Tours,
France}}

\date{}
\maketitle
\begin{center} {\bf\small Abstract}

\vspace{3mm}
\hspace{.05in}
\parbox{4.5in} {{\small Let $h:[0,\infty)\mapsto [0,\infty)$ be continuous and nondecreasing,
$h(t)>0$ if $t>0$, and $m,q$ be positive real numbers. We investigate the behavior when
$k\to\infty$ of the fundamental solutions $u=u_{k}$ of
$\prt_{t} u-\Gd u^m+h(t)u^q=0$ in $\Gw\ti (0,T)$ satisfying $u_{k}(x,0)=k\gd_0$.
The main question is wether the limit is still a solution of the above equation with an
isolated singularity at $(0,0)$, or a solution of the associated ordinary differential equation $ u'+h(t)u^q=0$ which
blows-up at $t=0$.
}}
\end{center}

\section {Introduction}
\setcounter{equation}
{0}

Let  $m$ and $q$ positive parameters and $h: [0,\infty)\mapsto [0,\infty)$ a nondecreasing
continuous. If one consider a reaction-diffusion equation such as
\begin {equation}
\label {intr1}
\prt_{t} u-\Gd u^m+h(t)u^q=0
\end {equation}
($u>0$ for simplicity) in a cylindrical domain $Q^T=\BBR^N\ti
(0,T)$ ($N\geq 1$), the behaviour of $u$ is subject to two
competing features: the diffusion associated to the partial
differential operator, here $-\Gd$, and the absorption which is
represented by the term $h(t)u^q$. When $q>1$ and $h(t)>0$ for
$t>0$, the absorption term is strong enough in order  positive
solution to satisfy an universal bound
\begin {equation}
\label {intr1'} 0\leq u(x,t)\leq U_h(t)=\left((q-1)\myint {0}{t}h(s)\,ds\right)^{-1/(q-1)}
\end {equation}
for every $(x,t)\in Q^T$. In addition, the function $U_h$ which
appears above is a particular solution of (\ref {intr1}). The
associated diffusion equation
\begin{equation}
\label {intr2}
\prt_{t} v-\Gd v^m=0
\end{equation}
admits fundamental solutions $v=v_k$ ($k>0$) which satisfy $v_k(x,0)=k\gd_0$ if
$m>(N-2)_+/N$. If
\begin{equation}
\label {intr3}
\myint{0}{T}\myint{B_{R}}{}h(t)v^q_k\,dx\,dt<\infty,\quad
B_R:=\{|x|<R\},
\end{equation}
for any $R\in (0,\infty]$, it is shown that (\ref {intr1}) admits
fundamental solutions $u=u_k$ in $Q^T$ which satisfy initial
condition $u_k(x,0)=k\gd_0$. The maximum principle holds and
therefore the mapping $k\mapsto u_k$ is increasing. If $h>0$ on
$(0,\infty)$ then due to universal bound\eqref{intr1'} there
exists $u_\infty=\lim_{k\to\infty}u_k$, and $u_\infty$ is a
solution of (\ref{intr1}) in $Q^{T}$. A natural question is
whether $u_\infty$ admits a singularity only at the origin $(0,0)$
or at other points too. Actually, in the last case it will imply
$u_{\infty}\equiv U$ since the following { alternative}
occurs:\medskip

\noindent (i) either $u_\infty=U$. ({\it complete initial
blow-up});\medskip

\noindent (ii) or $u_\infty$ is a solution singular at $(0,0)$ and
such that $\lim_{t\to 0}u(x,t)=0$ for all $x\neq 0$. ({\it
single-point initial blow-up}).\medskip

 This phenomenon is observed for the first
time by Marcus and V\'eron. They considered the semilinear
equation
\begin {equation}\label {heat1}
\prt_{t} u-\Gd u+h(t)u^q=0
\end {equation}
and  proved \cite [Prop. 5.2]{MV1}\medskip

\begin{theorem}\label{T1.1} {\it If $h(t)= e^{-\gk/t}$ ($\gk>0$), then the complete initial
blow-up occurs}.
\end{theorem}

However they raised the question whether this type of degeneracy
of the absorption is sharp or not. The method of \cite{MV1} relies
on the construction of subsolutions associated to very singular
solutions of equations
\begin {equation}
\label {intr4}
\prt_{t} u-\Gd u+c_\ge t^\alpha u^q=0
\end {equation}
for suitable $\ga>0$ and $c_\ge>0$, and on the study of
asymptotics of these solutions. One  the main result of present
paper states that if the degeneracy of the absorption terms is
lightly smaller respectivelly to Th. \ref{T1.1}, then localization
occurs.

\begin{theorem}\label{T1.2} {\it  If $h(t)=\exp (-\gw(t)/t)$, where $\gw$ is continuous,
nondecreasing and satisfies
\begin{equation}\label{Dini1} \int_0^1\myfrac {\sqrt\gw(s)}{s}ds<\infty,
\end {equation}
then $u_\infty$ has single-point initial blow-up at $(0,0)$}.
\end{theorem}

The method of the proof is totally different from the one of
Marcus and V\'eron  and based upon local energy estimates in the
spirit of the famous Saint-Venant 's principle (see
\cite{GSh1,OR,OI}). Using appropriate test functions we prove by
induction that the energy of the fundamental solutions  $u_{k}$
remains uniformly locally bounded in $\overline
Q^T\setminus\{(0,0)\}$.
\medskip

In the case of equation
\begin {equation}
\label {intr4}
\prt_{t} u-\Gd u+h(t)(e^{u}-1)=0
\end {equation}
the same type of phenomenon occurs, but at a different scale of degeneracy. We prove the following\medskip

\begin{theorem}\label{T1.3} {\it 1) If $h(t)=e^{-e^{\gk/t}}$ for some $\gk>0$, then the
complete initial blow-up occurs}.\smallskip

\noindent {\it 2)  If $h(t)=e^{-e^{\gw(t)/t}}$ for some $\gw\in
C(0,\infty)$ positive,  nondecreasing and satisfying
(\ref{Dini1}), then $u_\infty$ has single-point initial blow-up at
$(0,0)$}.\medskip
\end{theorem}

In this paper we also extend the study of equation (\ref{intr1})
to the case $m\neq 1$. The situation differs completely
corresponding to $m>1$, the porous media equation with slow
diffusion, and to $(N-2)_+/N<m<1$, the fast diffusion equation.
Concerning the porous media equation, we prove \medskip

\begin{theorem}\label{T1.4}{\it If $q>m>1$ and $h$ is nondecreasing and satisfies
$h(t)=O(t^{(q-m)/(m-1)})$ as $t\to 0$, then $u_{\ity}\equiv
U_{h}$.}
\end{theorem}

We give two proofs. The first one, valid only in the subscritical
case $1<m<q<m+2/N$, is based upon the construction of suitable
subsolutions, as in the semilinear case. The second one, based
upon scaling transformations, is valid in all the cases $q+1>2m>2$
where the $u_{k}$ exists. It reduces to proving that the equation
$$-\Gd \Psi-\Psi^{1/m}+\Psi^{q/m}=0\quad\mbox {in }\BBR^N
$$
admits only one positive solution, the constant $1$.
The localization counter part is as follows,\medskip

\begin{theorem}\label{T1.5}{\it Assume $q>m>1$, in Equation (\ref{intr1}). If
$h(t)=t^{(q-m)/(m-1)}\gw^{-1}(t)$ with $\gw(t)\to 0$ as $t\to 0$,
and
\begin{equation}
\label{intr6*} \int_0^1\gw^\gth(s)\myfrac {ds}{s}<\infty
\end{equation}
where
$$
\gth=\myfrac{m^2-1}{[N(m-1)+2(m+1)](q-1)},
$$ then
$u_\infty$ has single-point initial blow-up at $0,0)$.}
\end{theorem}

Actually, the method is applicable to a much more general class of
equations.\medskip

In the fast diffusion case there is always localization.\medskip

\begin{theorem}\label{T1.6}{\it  Assume $(N-2)_+/N<m<1$ and $q>1$, in Equation
(\ref{intr1}). Then
\begin {equation}
\label {intr6}
u_{\infty}(x,t)\leq \min\left\{U_{h}(t),C_{*}\left(\myfrac {t}{\abs x^2}\right)
^{1/1-m)}\right\}
\end {equation}
where
$$C_{*}=\left(\myfrac{(1-m)^3}{2m(mN+2-N}\right)^{1/(1-m)}.
$$
}\end{theorem}

This type of problem has an elliptic counterpart which is initiated in \cite {MV3} where the following question is considered: suppose $\Gw$ is a $C^2$ bounded domain in $\BBR^N$, $q>1$ and $h\in C(0,\infty)$ is positive. What is the limit, when $k\to\infty$ of the solutions (when they exist) $u=u_{k}$ of the following problem
\begin {equation}\label {elliptic}\left\{\BA {l}
-\Gd u+h(\gr(x))u^q=0\quad\mbox {in }\Gw\\[2mm]
\phantom{-\Gd u+h(\gr(x))u} u=k\gd_{0}\quad\mbox {in }\prt\Gw,
\EA\right.\end {equation} where $\gr(x)=\dist (x,\prt\Gw)$. It is
proved in \cite {MV3} that, if $h(t)=e^{-1/t}$, then
$u_{\infty}$($:=\lim_{k\to\infty}u_{k}$) is the maximal solution
of the equation in $\Gw$, that is the function which satisfies
\begin {equation}\label {elliptic2}\left\{\BA {l}
-\Gd u+h(\gr(x))u^q=0\quad\mbox {in }\Gw\\[2mm]
\phantom{,-} \lim_{\gr(x)\to 0}u(x)=\infty.
\EA\right.
\end{equation}
On the contrary, if $h(t)=t^\ga$, for $\ga>0$ and
$1<q<(N+1+\ga)/(N-1)$, it is proved in \cite{MV4} that
$u_{\infty}$ has an isolated singularity at $0$, and vanishes
everywhere outside $0$. In a forthcoming article we shall study
this localization of singularity
phenomenon for
the complete nonlinear elliptic problem,
replacing the powers by more general functions, and the ordinary Laplacian
by the $p$-Laplacian operator.\\

\noindent Our paper is organized as follows: \S1 Introduction. In
\S 2 we study sufficient conditions of complete initial blow-up
for semilinear heat equation. In \S 3 we prove sharp sufficient
condition of existence of single point initial blow-up for heat
equation with power nonlinear absorption. In \S 4 local energy
method from \S 3 is adapted to the heat equation with nonpower
absorption nonlinearity. \S 5 deals with porous media equation
with power nonlinear absorption, \S 6 --- the fast diffusion
equation with nonlinear absorption.

\section{Complete initial blow-up for semilinear\\ heat equation}
\setcounter{equation}{0}

We recall the standard result concerning the existence of a fundamental solution $u=u_{k}$ ($k>0$) to the following problem
\begin{equation}
\label{gen}\left\{\BA {l}
\partial_{t}u-\Gd u+g(x,t,u)=0\quad \mbox {in }Q^{T}=\BBR^N\ti (0,T)\\[2mm]
u(x,0)=k\gd_{0}. \EA\right.
\end{equation}
If $v$ is defined in
$Q^{T}$, we denote by $\tilde g(v)$ the function $(x,t)\mapsto
g(x,t,v(x,t))$. By a solution we mean a function $u\in
L^1_{loc}(\overline Q^{T})$ such that $\tilde g(u)\in
L^1_{loc}(\overline Q^{T})$, which verifies
\begin{equation}
\label{weak1} \dint_{Q^{T}}\left(-u\prt_{t}\gf-u\Gd\gf+\tilde
g(u)\gf\right)dxdt=k\gf(0,0),
\end{equation}
for any $\gf\in C_{0}^{2,1}(\BBR^N\ti [0,T)\ti\BBR)$. We denote by
$E(x,t)=(4\gp t)^{-N/2}e^{-\abs x^2/4t}$ the fundamental solution
of the heat equation in $Q^\infty$, by $B_R(a)$ an open ball of
center $a$ and radius $R$, and $B_R(0)=B_R$. The following result
is classical

\bth{exsth} Let $g\in C(\BBR^N\ti [0,T]\ti\BBR)$ such that
$g(x,t,r)\geq 0$ on $\BBR^N\ti [0,T]\ti\BBR_{+}$, and assume that
$g=g_{1}+g_{2}$ where $g_{1}$ and $g_{2}$ are respectively
nondecreasing and locally Lipschitz continuous with respect to the
$r$-variable functions. Let $k>0$ be such that
\begin {equation}
\label {gen-k}
\myint{0}{T}\myint{B_{R}}{}g(x,t,kE(x,t))dxdt<\infty.
\end {equation}
for any $R>0$. Then there exists a solution $u=u_{k}$ to problem
\eqref{gen}. Furthermore, if $g_{2}=0$, then $u_{k}$ is unique.
\es

Function $g(x,t,r)=e^{-\kappa/t} {\abs r}^{q-1}r$, with $\kappa
>0$ and $q>1$,  satisfies \eqref{gen-k}. Thus the problem
\begin {equation}
\label {expu^q}\left\{\BA {l}
\partial_{t}u-\Gd u+e^{-\kappa/t}\abs u^{q-1}u=0\quad \mbox {in }Q^{\infty}\\[2mm]
u(x,0)=k\gd_{0}.
\EA\right.\end {equation}
admits a unique solution. The next result is proved in \cite {MV1}, but we recall the proof both  for the sake of completeness and to present the key-lines of the method in a simple case.

\bth{powerth}\label{initial blow-up} For $k>0$, let $u_{k}$ denote
the solution of (\ref {expu^q}) in $Q^{\infty}$. Then
$u_{k}\uparrow U_{S}$ as $k\to\infty$, where
\begin {equation}
\label {expu^q 2}
U_{S}(t)=\left((q-1)\int_{0}^te^{-\kappa/s}ds\right)^{1/(1-q)}\,,\forevery t>0.
\end {equation}
\es
\Proof {\it Case 1.} $1<q<1+2/N$. For any $\epsilon>0$, $u_{k}=u$ satisfies
\begin {equation}
\label {expu^q epsilon}
\partial_{t}u-\Gd u+e^{-\kappa/\epsilon}u^q\geq 0
\end {equation}
on $Q^{\ge}$. Therefore if $v=v_{k}$ is the solution of
\begin {equation}
\label {expv^q epsilon}\left\{\BA {l}
\partial_{t}v-\Gd v+e^{-\kappa/\epsilon}v^q= 0\quad \mbox {in }Q^{\infty}\\[2mm]
v(x,0)=k\gd_{0}, \EA\right.\end {equation} there holds $u_{k}\geq
v_{k}$. Passage to the limit $k\to\infty$,  yields
\begin {equation}
\label {u>v 2} \lim_{k\to\infty}u_{k}:=u_{\infty}\geq v_{\infty}=\lim_{k\to\infty}v_{k}\;\mbox { in } Q^{\ge}.
\end {equation}
If we write $v_{\infty}(x,t)=e^{\kappa/\epsilon(q-1)}t^{-1/(q-1)}f(x/\sqrt
t)$, then $f$ is radial and  satisfies
$$\left\{\BA {l}
f''+\left(\myfrac {N-1}{r}+\myfrac {r}{2}\right)f'+\myfrac
{1}{q-1}f-f^q=0\quad\mbox { on } (0,\infty),\\
f'(0)=0\,,\;\lim_{r\to\infty}r^{2/q-1)}f(r)=0. \EA\right.$$
Furthermore the asymptotics of $f$ is given in \cite {BPT},
$$
f(r)=Cr^{2/(q-1)-N}e^{-r^{2}/4}(1+\circ (1)))\,,\;\mbox { as } r\to\infty, $$
for some $C=C(N,q)>0$. Therefore
\begin {equation}
\label {f} f(r)\geq \tilde C(r+1)^{2/(q-1)-N}e^{-r^{2}/4}\,\forevery r\geq 0,
\end {equation}
for some $\tilde C=\tilde C(N,q)>0$. If we take $t=\epsilon$, we derive from (\ref {u>v 2})
\begin {equation}
\label {u>v 3} u_{\infty}(x,t)\geq e^{\kappa/t(q-1)}t^{-1/(q-1)}f(x/\sqrt t)\;\mbox { in
} \mathbb R^{N}.
\end {equation}
Let $0<\ell<2\sqrt{\kappa/(q-1)}$. Inequalities (\ref {f}) and (\ref {u>v 3}) imply
\begin {equation}
\label {u>v 4} u_{\infty}(x,t)\geq \tilde
Ct^{-1/(q-1)}e^{(\kappa/(q-1)-\ell^{2}/4)t^{-1}}\,,\forevery x\in \bar B_{\ell}.
\end {equation}
Therefore $\lim_{t\to 0}u_{\infty}(x,t)=\infty\,,\forevery x\in \bar B_{\ell}$. We
 pick some point $x_{0}$ in $B_{\ell}$. Since for any $k>0$, the solution
$u_{k\gd_{x_{0}}}$ of (\ref {expu^q}) with initial value $k\gd_{x_{0}}$ can be
approximated by solutions with bounded initial data  and support in
$B_{\sigma}(x_{0})$ ($0<\gs<\ell-\abs {x_{0}}$), the previous inequality implies $$
u_{\infty}(x,t)\geq u_{\infty}(x-x_{0},t). $$
Reversing the role of $0$ and $x_{0}$ yields to $$
u_{\infty}(x,t)=u_{\infty}(x-x_{0},t). $$
If we iterate this process we derive
\begin {equation}
\label {initial expl} u_{\infty}(x,t)=u_{\infty}(x-y,t)\,,\forevery y\in \mathbb
R^{N}.
\end {equation}
\medskip Since $u_{k\delta_{y}}$ is radial with respect to $y$,
(\ref {initial expl}) implies that $u_{\infty}(x,t)$ is
independent of $x$ and therefore it is solution of
\begin {equation}
\label {ODE} \left\{\BA {l}
z'+e^{-\kappa/t}z^q=0\quad\mbox { on } (0,\infty)\\[2mm]
\,\,\,\lim_{t\to 0}z(t)=\infty.
\EA\right.
\end {equation}
Thus $u_{\infty}=U_{S}$ where $U_{S}$ is defined by (\ref {expu^q 2}).\medskip

\noindent {\it Case 2.} $q\geq 1+2/N$. Let $\ga>0$ such that $q<q_{c,\ga}=1+2(1+\ga)/N$.
We write $e^{-\kappa/t}=t^{\ga}\tilde h(t)$
with $
\tilde h(t)=t^{-\ga}e^{-\kappa/t}$.
The function $\tilde h$ is increasing on $(0,\kappa/\ga]$ and we extend it by $\tilde
h(0)=0$. Let $0<\epsilon\leq \kappa/\ga$, then the solution $u=u_{k}$ of (\ref {expu^q}) verifies
$$
\partial_{t}u-\Gd u+\tilde h(\epsilon)t^\ga u^q\geq 0, $$
in $\mathbb R^{N}\times (0,\epsilon]$. As in Case 1, $u$ is bounded from below on $\mathbb
R^{N}\times (0,\epsilon]$ by $\left(\tilde h(\epsilon)\right)^{-1/(q-1)}v_{\infty}$ where
$v_{\infty}=v$ is
is the very singular solution of
\begin {equation}\label {VSS-a}
\prt _{t}v-\Gd v+t^\ga v^q=0.
\end {equation}
Then $v_{\infty}(x,t)=t^{-(1+\ga)/(q-1)}f_{\ga}(\abs x/\sqrt t)$, and $f_{\ga}=f$ satisfies
$$\left\{\BA {l}
f''+\left(\myfrac {N-1}{r}+\myfrac {r}{2}\right)f'+\myfrac
{1+\ga}{q-1}f-f^q=0\quad\mbox { on } (0,\infty),\\
f'(0)=0\,,\;\lim_{r\to\infty}r^{2(1+\ga)/q-1)}f(r)=0. \EA\right.$$
The asymptotics of $f_{\ga}$ is given in \cite {MV2}
$$f_{\ga}(r)=Cr^{2(1+\ga)/(q-1)-N}e^{-r^2/4}(1+\circ (1))\quad\mbox {as }r\to\infty,
$$
thus
$$f_{\ga}(r)\geq \tilde C(1+r)^{2(1+\ga)/(q-1)-N}e^{-r^2/4}\forevery r\in \BBR_{+}.
$$
Consequently
\begin {equation}
\label {u>v 41} u(x,t)\geq \tilde C
e^{(\kappa/(q-1)-\ell^{2}/4)t^{-1}}\,,\forevery x\in \bar B_{\ell}.
\end {equation}
Taking again $0<\ell<2\sqrt{\kappa/(q-1)}$, we derive
$$
\lim_{t\to 0}u(x,t)=\infty\,,\forevery x\in \bar B_{\ell}. $$ As
in the Case 1, it yields to $u_{\infty}(x,t)=u_{\infty}(x-y,t)$
for any $y\in\mathbb R^{N}$, and finally
$u_{\infty}(x,t)=U_{S}(t)$.\qeda \medskip

Next we consider Cauchy problem for diffusion equation with an
exponential type absorption term
\begin {equation}\label{exp1}\left\{\BA {l}
\prt _{t}u-\Gd u+h(t)e^{u}=0\quad\mbox {in }Q^{\infty}\\[2mm]
u(x,0)=k\gd_{0}
\EA\right.
\end {equation}
 where $h\in C(\BBR_{+})$ is nonnegative.  \rth{exsth} yields the following existence
 result:
\bprop {exp-exist}Assume $h$ satisfies
\begin {equation}\label{exp2}
\lim_{t\to 0}t^{N/2}\ln h(t)=-\infty.
\end {equation}
Then for any $k>0$ problem (\ref {exp1}) admits a unique solution $u=u_{k}$. Furthermore
\begin {equation}\label{exp3}
u_{k}(x,t)\leq V_{S}(t):=-\ln\left(\myint{0}{t}h(s)ds\right)\forevery (x,t)\in Q^{\infty}.
\end {equation}
\es

Notice that estimate (\ref{exp3}) is a consequence of the fact
that $V_{S}$ satisfies the associated O.D.E.
$$y'+h(t)e^{y}=0\quad\mbox {in }(0,\infty),
$$
with infinite initial value. Our main result concerning
nonexistence of localized singularities for equation \eqref{exp1}
is
\bth {exp-th} Let $h(t)=e^{-e^{\gs/t}}$ for some $\gs>0$ and any $t>0$. Then $u_{k}\uparrow V_{S}$ as $k\to\infty$.
\es
\Proof {\it Step 1. Construction of an approximate very singular solution. }
For $n>1$ and $c_{n}>0$ to be defined later on, let $v=V_n$ be the very singular solution  of
\begin {equation}\label{main-k}
\prt _tv-\Gd v+c_nt^{\ga_n}v^n=0.
\end {equation}
The necessary and sufficient condition for the existence of a $V_{n}$     is
$$n<1+N(\ga_n+1)/2.
$$
This function is obtained in the form
$$V_{n}(x,t)=t^{-(1+\ga_n)/(n-1)}F(x/\sqrt t),
$$
where $F$ solves
$$\Gd F+\myfrac {1}{2}\xi.DF+\myfrac {1+\ga_n}{n-1}F-c_nF^n=0.
$$
We fix
\begin {equation}\la{fix-n}
\myfrac {1+\ga_n}{n-1}=1+\myfrac {N}{2}\Longleftrightarrow
\ga_n=(2+N)(n-1)/2-1,
\end {equation}
and set
$$f_n=c_n^{1/(n-1)}F.
$$
Then $f_n$ solves
$$\Gd f_n+\myfrac {1}{2}\xi.Df_n+\myfrac {N+2}{2}f_n-f_n^n=0.
$$
We prove that $f_n$ has an asymptotic expansion essentially
independent of $n$, in the following form
\begin {equation}\la{asymp}
f_n(\xi)\geq \gd (\abs \xi^2+1)e^{-\abs{\xi}^{2}/4}
\Longrightarrow V_n(x,t)\geq  \gd c_n^{-1/(n-1)} t^{-2-N/2}(\abs
x^2+t)e^{-\abs x^2/4t}
\end{equation}
It order to see that, we put
$$\tilde f_{n}= \left(\myfrac {2}{N+2}\right)^{1/(n-1)}f_{n}
$$
then
$$\Gd \tilde f_n+\myfrac {1}{2}\xi.D\tilde f_n+\myfrac {N+2}{2}\tilde f_n-\myfrac {N+2}{2}\tilde f_n^n=0.
$$
By the maximum principle $0\leq \tilde f_{n}\leq 1$ so that $0\leq \tilde f^{n'}_{n}\leq \tilde f^n_{n}$ for $n'>n$. Thus
$$\Gd \tilde f_n+\myfrac {1}{2}\xi.D\tilde f_n+\myfrac {N+2}{2}\tilde f_n-\myfrac {N+2}{2}\tilde f_{n}^{n'}\geq 0,
$$
which implies that $\tilde f_{n}$ is a subsolution of the equation
for $\tilde f_{n'}$ and therefore,
\begin {equation}\la{asymp1}
n'>n\Longrightarrow \tilde f_{n}\leq \tilde f_{n'}\Longleftrightarrow f_{n}\leq
\left(\myfrac {N+2}{2}\right)^{(n'-n)/(n-1)(n'-1)} f_{n'}.
\end{equation}
In the particular case $n=n^*=(N+4)/(N+2)$, the equation falls
into the scoop of Brezis-Peletier-Terman study since it can also
be written in the form
$$\Gd f_{n^*}+\myfrac {1}{2}\xi.Df_{n^*}+\myfrac {1}{n^*-1}f_{n^*}-f_{n^*}^{n^*}=0.
$$
and their asymptotic expansion applies  (with $2/(n^*-1)-N=2$) as $\abs \xi \to\infty$:
\begin {equation}\la{asymp2}
f_{n^*}(\xi)=C\abs \xi^2e^{-\abs \xi^2/4}(1+\circ (1))\Longrightarrow f_{n^*}(\xi)\geq\gd_{*}(\abs \xi^2+1)e^{-\abs \xi^2/4}\forevery \xi.
\end{equation}
Combining (\ref{asymp1}) with $n=n^*$ and $n'$ replaced by $n$, and (\ref{asymp2}), we get
\begin {equation}\la{asymp3}
 f_{n}(\xi)\geq\gd_{*}\left(\myfrac {2}{N+2}\right)^{(n-n^*)/(n-1)(n^*-1)}(\abs \xi^2+1)e^{-\abs \xi^2/4}\forevery \xi.
\end{equation}
Since $n\mapsto\left(2/(N+2\right)^{(n-n^*)/(n-1)(n^*-1)}$ is bounded from below independently of $n>n^*$, we get
(\ref{asymp}).\medskip

\noindent{\it Step 2. Some estimates from below for a related problem. }
In order to have $v_n\leq u$ in the range of value of $u$, which is
\begin {equation}\la{sup-V}
u(t)\leq V_{S}(t)=-\ln\left(\myint{0}{t}h(s)ds\right)\forevery t>0,
\end{equation}
we need $v=v_{n}$ to be a subsolution near $t=0$ of the equation that $u$ verifies. Furthermore this can be done up to some bounded function. It is sufficient to have
\begin {equation} \label {ineq'}
c_nt^{\ga_n}(x^n+1)\geq h(t)e^x,\forevery t\in (0,\gt_n],\;x\in [0,V_{S}(t)]
\end {equation}
where $\gt_n$ has to be defined. In particular, at the end points of the interval,
\begin{equation} \label{ineq0}\begin{cases}
(i)\;c_nt^{\ga_k}\geq h(t)\\
(ii)\;
c_nt^{\ga_n}\left(\ln^n\left(\dfrac{1}{\int_0^ta(s)ds}\right)+1\right)
\geq \dfrac{h(t)}{\int_0^t h(s)\,ds}. \end{cases}
\end{equation}
 We  write
(\ref{ineq'}) in the form
\begin {equation} \label {ineq}
\myfrac {e^x}{1+x^n}\leq\myfrac {c_nt^{\ga_n}}{h(t)},
\end {equation}
and set
$$\gf(x)=\myfrac {e^x}{1+x^n}.$$
Then
$$\gf'(x)=e^x\myfrac {1+x^n-nx^{n-1}}{(1+x^n)^2}.
$$
The sign of $\gf'$ is the same as the one of
$\psi(x)=1+x^k-nx^{n-1},
$
a function which decreasing then increasing, is positive near $0$, vanishes somewhere between $0$ and $1$ and again between $n-1$ and $n$. The first maximum of $\gf$ is less than $e/2$. This is not important in (\ref {ineq}) since we can always assume that the minimum of $c_kt^{\ga_k}/h(t)$ is larger than $e/2$. Therefore, it is sufficient to have
\begin {equation} \label {ineq2}
\myfrac {e^{V_{S}(t)}}{1+V_{S}^n(t)}\leq\myfrac {c_nt^{\ga_n}}{h(t)},
\end {equation}
in order to have (\ref {ineq}). This is exactly (\ref {ineq0})-ii.
If we express $h(t)$ in the form
$$h(t)=-\gw'(t)e^{-\gw(t)},$$
then (\ref {ineq0})-ii is equivalent to
\begin {equation} \label {ineq3}
c_nt^{\ga_n}\left(\gw^n(t)+1\right)
\geq -\gw'(t).
\end {equation}
Since
$$\gw^n(t)+1\geq 2^{1-n}(\gw(t)+1)^n,
$$
we  associate the following O. D. E. on $\BBR_{+}$
$$c_nt^{\ga_n}=2^{1-n}\myfrac {-\eta'}{(\eta+1)^n},
$$
the maximal solution of which is
$$\eta (t)=\myfrac {1}{2}\left(\myfrac {1}{c_n(n-1)}\right)^{1/(n-1)}
t^{-(\ga_n+1)/(n-1)}=\myfrac {1}{2}\left(\myfrac {1}{c_n(n-1)}\right)^{1/(n-1)}
t^{-1-N/2}.
$$
If we write $\gw$ in the form
$$\gw(t)=e^{\ga(t)},
$$
with $\ga(0)=\infty$, $\ga'<0$, then (\ref {ineq0})-ii becomes
$$c_nt^{\ga_n}\left(e^{n \ga(t)}+1\right)\geq-\ga'(t)e^{\ga(t)},
$$
and this inequality is ensured provided
\begin {equation}\label {choi0}
c_nt^{\ga_n}e^{(n-1) \ga(t)}\geq-\ga'(t)\Longleftrightarrow c_{n}\geq -\ga'(t)e^{(1-n)\ga(t)-\ga_{n}\ln t}=
 -t\ga'(t)e^{(1-n)\left(\ga(t)+2^{-1}(N+2)\ln t\right)},
\end {equation}
by replacing $\ga_{n}$ by its value. Next we fix
\begin {equation}\label {choi}
\ga(t)=\ga_{\gs}(t)=\myfrac {\gs}{t}\forevery t>0
\end {equation}
where $\gs>0$ is a parameter, thus
$$
 -t\ga'(t)e^{(1-n)\left(\ga(t)+2^{-1}(N+2)\ln t\right)}=
 e^{(1-n)\gs/t-\left(2^{-1}(n-1)(N+2)+1\right)\ln t}=e^{\gr(t)}.
$$
In order to have (\ref {choi0}) it is sufficient to have the monotonicity of the function $\gr$ and
$$\gr'(t)=\myfrac {\gs(n-1)}{t^2}-\myfrac {n(N+2)-N}{2t}
$$
Then there exist $\gga>0$, independent of $k$ and $\gs$ such that $\gr'(t)>0$ on $(0,\gs\gga]$. Consequently, inequality (\ref {choi0}) is ensured on $(0,\ge]\subset (0,\gs\gga]$ as soon as
\begin {equation}\label {choi1}
c_{n}\geq e^{\gr(\ge)}=e^{(1-n)\gs/\ge-2^{-1}(n(N+2)-N)\ln\ge}.
\end {equation}
{\it Step 3. Complete initial blow-up for a related problem.}
Assume now
\begin {equation}\label {choi1}
h(t)=\tilde\gs t^{-2}e^{\tilde\gs t^{-1}-e^{\tilde\gs/t}}
\end {equation}
 for some $\tilde\gs>0$. For $n>2$, we fix $\ge<\tilde\gs\gga$ and take $c_{n}=e^{\gr(\ge)}$. On $(0,\ge]$ we have
$$c_{n}t^{\ga_{n}}(e^{n\ga(t)}+1)\geq -\ga'(t)e^{\ga(t)}.
$$
Therefore, if $u=u_{k}$ is  the solution of (\ref {exp1}) with $h(t)$ given by (\ref{choi1}), it satisfies $u(t)\leq V_{S}(t)$, where $V_{S}$ is given by (\ref{sup-V}), and
$$ \prt_{t}u-\Gd u+c_{n}t^{\ga_{n}}(u^n+1)\geq 0\quad\mbox {in }Q^{\ge}.
$$
Therefore $u$ is larger that the solution $v=\tilde v_{k}$ of
$$ \prt_{t}v-\Gd v+c_{n}t^{\ga_{n}}(v^n+1)= 0\quad\mbox {in }Q^{\ge},
$$
with $\tilde v_{k}(0)=k\gd_{0}$. Furthermore $\tilde v_{k}\geq  v_{k}-c_{n}t^{\ga_{n}+1}/(\ga_{n}+1)$, where
$v= v_{k}$ solves
$$ \prt_{t}v-\Gd v+c_{n}t^{\ga_{n}}v^n= 0\quad\mbox {in }Q^{\ge},
$$
with $ v_{k}(0)=k\gd_{0}$. If we let $k\to\infty$, we derive from (\ref{asymp}) and by replacing $c_{n}=e^{\gr(\ge)}$ by
its precise value $e^{(1-n)\gs/\ge-2^{-1}(n(N+2)-N)\ln\ge}$, that
$$u_{\infty}(x,t)\geq V_{n}(x,t)-\myfrac {c_{n}t^{\ga_{n}+1}}{\ga_{n}+1}
\geq \gd t^{-2-N/2}(\abs x^2+t)e^{\frac {\gs}{\ge}+\frac {(n(N+2)-N\ln\ge}{n-1}-\frac {\abs x^2}{4t}}
$$
on $(0,\ge]$. In particular
\begin {equation}\label {choi2}
u_{\infty}(x,\ge)\geq
\gd \ge^{-2-N/2}(\abs x^2+\ge)e^{\frac {\gs}{\ge}+\frac {(n(N+2)-N\ln\ge}{n-1}-\frac {\abs x^2}{4\ge}}.
\end {equation}
Taking $\abs x^2<\gs/4$ yields to
$$\lim_{\ge\to 0}\ge^{-2-N/2}(\abs x^2+\ge)e^{\frac {\gs}{\ge}+\frac {(n(N+2)-N\ln\ge}{n-1}-\frac {\abs x^2}{4\ge}}=\infty.
$$
Thus
$$\lim_{\ge \to 0}u_{\infty}(x,\ge)=\infty,\forevery x\in B_{\sqrt\gs/2}.
$$
As in the proof of \rth {powerth}, it implies $u_{\infty}=V_{S}$. \medskip

\noindent{\it Step 4. End of the proof. } Since for any $\gs>\tilde \gs>0$ there exists an interval $(0,\gth]$ on which
$$\tilde \gs t^{-2}e^{\gs' t^{-1}-e^{\gs'/t}}\geq e^{-e^{\gs/t}},$$
any solution of (\ref {exp1}) with $h(t)$ given by (\ref{choi1}) is a subsolution in $Q^{\gth}$ of the same equation with $h(t)=e^{-e^{-\gs/t}}$. This implies the claim.\qeda
\section{Single point initial blow-up for semilinear\\ heat equation}
\setcounter{equation}{0}

We  consider the following Cauchy problem
\begin {equation}
\label {au^q}\left\{\BA {l}
\partial_{t}u-\Gd u+h(t)\abs u^{q-1}u=0\quad \mbox {in }Q^{\infty}\\[2mm]
u(x,0)=k\gd_{0}. \EA\right.\end {equation} The first result
dealing with the localization of the blow-up that we prove is the
following. \bth {locth1} Assume $h(t)=e^{-\gw(t)/t}$ where $\gw\in
C([0,\infty))$ is positive, nondecreasing  function which
satisfies $\gw(s)\geq s^{\ga_{0}}$ for some $\ga_{0}\in [0,1)$ and
any $s>0$, and the following Dini like condition holds:
\begin {equation}
\label {dini1}
\myint{0}{1}\myfrac{\sqrt{\gw(s)}}{s}ds<\infty.
\end {equation}
Then $u_{k}$ always exists and
$u_{\infty}:=\lim_{k\to\infty}u_{k}$ has a point-wise singularity
at $(0,0)$. \es \Proof The proof is based on the study of
asymptotic properties as $k\to \infty$ of solutions $u= u_k$ of
the regularized Cauchy problem
\begin{equation}\label{2.1}\left\{\BA {l}
u_t-\Gd u+h(t)|u|^{q-1}u=0\quad\mbox{in } Q^T,\\[2mm]
u(x,0)=u_{0,k}(x)=M^{1/2}_kk^{- N/2}\gd_k(x)\forevery x\in\BBR^N,
\EA\right.\end{equation} where $\gd_k\in C(\BBR^N),\ \supp
\gd_k\subset\curlybr{|x|\le k^{-1}},\ \gd_k\warrow\gd(x)$ weakly
in the sense of measures as $k\to \infty$ and $\curlybr{M_k}$ is
some sequence tending to $\infty$ as $k\to \infty$ fast enough so
that
\begin{equation}\label{2.1.2}
M_k^{1/2}k^{-N/2}\to \infty \text{ as } k\to \infty.
\end{equation}
Without loss of generality we will suppose that
\begin{equation}\label{2.3}
\| \gd_k(x) \|_{L_2(\BBR^N)}^2\le c_0k^N\quad \forall\,
k\in\BBN,\quad c_0=\text{const}.
\end{equation}

Our method of analysis is some variant of the local energy
estimates method (also called Saint-Venant principle), developed,
particulary, in \cite{OI,OR,SS,Sh1,Sh2} (see also review in
\cite{GSh1}). Let introduce the families of subdomains
\begin{align*}
& \Gw(\gt)=\BBR^N\cap\curlybr{|x|>\gt}\quad \forall\, \gt>0,\\
& Q^r(\gt)=\Gw(\gt)\times(0,r)\quad \forall\, r\in(0,T),\\
& Q_r(\gt)=\Gw(\gt)\times(r,T) \quad \forall\,r\in(0,T).
\end{align*}
{\it Step 1. The local energy framework.} We fix arbitrary $k\in
\BBN$ and consider solution $u= u_k$ of \eqref{2.1}, but for
convenience we will denote it by $u$. Firstly we deduce some
integral vanishing properties of solution $u$ in the family of
subdomains $Q_r:=\mathbb{R}^N\times(r,T)$. Multiplying (\ref{2.1})
by $u(x,t)\exp\paran{-\dfrac{t-r}{1+T-r}}$ and integrating in
$Q_r$, we get
\begin{multline}\label{2.4}
\paran{2\exp\paran{\frac{T-r}{1+T-r}}}^{-1}\int_{\BBR^N}\abs{u(x,T)}^2dx\\
+\int_{Q_r}\paran{|D_xu|^2+h(t)|u|^{q+1}}\exp\paran{-\frac{t-r}{1+T-r}}dxdt\\
+\frac{1}{1+T-r}\int_{Q_r}|u|^2\exp\paran{-\frac{t-r}{1+T-r}}dxdt\\
=2^{-1}\int_{\Gw(\gt)}|u(x,r)|^2\,dx+2^{-1}\int_{\BBR^N\sms\Gw(\gt)}|u(x,r)|^2dx,
\end{multline}
where $\gt>0$ is arbitrary parameter.  Using
H\"older's inequality, it is easy to check  that
\begin{equation}\label{2.5}
\int_{\BBR^N\sms\Gw(\gt)}|u(x,r)|^2\,dx\le
c\gt^{\frac{N(q-1)}{q+1}}h(r)^{-\frac{2}{q+1}}\paran{\int_{\BBR^N\sms\Gw(\gt)}
|u(x,r)|^{q+1}h(r)\,dx}^{\frac{2}{q+1}}.
\end{equation}
Here and further we will denote by $c,\ c_i$ different positive
constants which do not depend on parameters $k,\ \gt,\ r$, but the precise value of which may change from one ocurrence to another.
Let us consider now the energy functions
\begin{equation}\label{2.5*}
 I_1(r)=\int_{Q_r}|D_xu|^2\,dx\,dt,\quad
I_2(r)=\int_{Q_r}h(t)|u(x,t)|^{q+1}\,dxdt,\quad
I_3(r)=\int_{Q_r}|u|^2\,dxdt.
\end{equation}
It is easy to check that
$$
-\frac{dI_2(r)}{dr}=\int_{\BBR^N}h(r)|u(x,r)|^{q+1}\,dx
\geq\int_{\BBR^N\sms\Gw(\gt)}h(r)|u(x,r)|^{q+1}\,dx\quad
\forall\,\gt>0.
$$
Therefore it follows from \eqref{2.4} and \eqref{2.5}
\begin{multline}\label{2.6}
\int_{\BBR^N}|u(x,T)|^2\,dx+I_1(r)+I_2(r)+I_3(r)\leq
c\gt^{\frac{N(q-1)}{q+1}}h(r)^{-\frac{2}{q+1}}\paran{-I'_2(r)}^{\frac2{q+1}}+
c\int_{\Gw(\gt)}|u(x,r)|^2\,dx\\ \forall\, \gt>0,\
\forall\,r:0<r<T.
\end{multline}
Next we introduce additional energy functions
\begin{equation}\label{2.6*}
f(r,\gt)=\int_{\Gw(\gt)}|u(x,r)|^2\,dx,\quad
E_1(r,\gt)=\int_{Q^r(\gt)}|D_xu|^2\,dxdt,\quad
E_2(r,\gt)=\int_{Q^r(\gt)}|u|^2\,dxdt.
\end{equation}
Now we deduce some vanishing estimates of these energy functions.
Let $\gm$ be some nondecreasing smooth function defined on
$(0,\infty)$, $\mu(\gt)>0$ for $\gt>0$ (a more precise definition
will be fixed later on). Then multiplying the equation \eqref{2.1}
by $u(x,t)\exp(-\gm^2(\gt)t)$ and integrating in domain $Q^r(\gt)$
with $\gt>k^{-1}$ (remember that $\supp
u_{0,k}\sbs\curlybr{|x|<k^{-1}}$) we deduce easily
\begin{multline}\label{2.7}
2^{-1}f_{\mu,r}(\gt)+J_{\gm,r}(\gt)
:=2^{-1}\int_{\Gw(\gt)}|u(x,r)|^2\exp(-\gm^2(\gt)r)\,dx+\\
\int_{Q^r(\gt)}\paran{|\nabla_xu|^2+\gm^2(\gt)|u|^2}\exp(-\gm^2(\gt)t)\,dxdt\\
\leq\mu(\gt)^{-1}\int_{\partial\Gw(\gt)\times(0,r)}\paran{|\nabla_xu|^2+\mu^2(\gt)|u|^2}
\exp(-\gm^2(\gt)t)\,dsdt\quad \forall\, \gt>k^{-1}.
\end{multline} Clearly there holds
\begin{multline*}
\frac{dJ_{\mu,r}(\gt)}{d\gt}=-\int_{\partial\Gw(\gt)\times(0,r)}
\paran{|\nabla_xu|^2+\mu^2(\gt)|u|^2}\exp(-\gm^2(\gt)t)\,dsdt\\+
\int_{Q^r(\gt)}2\mu\mu'(\gt)|u|^2\exp(-\gm^2(\gt)t)\,dxdt\\
-2\int_{Q^r(\gt)}\mu\mu'(\gt)t\paran{|\nabla_xu|^2
+\mu^2(\gt)|u|^2}\exp(-\gm^2(\gt)t)\,dxdt.
\end{multline*}
Since $\mu'(\gt)>0$, it follows from \eqref{2.7},
\begin{equation}\label{2.8}
2^{-1}f_{\mu,r}(\gt)+J_{\mu,r}(\gt)\leq\mu(\gt)^{-1}\sqbr{-\frac{d}{d\gt}J_{\mu,r}(\gt)+
2\int_{Q^r(\gt)}\mu(\gt)\mu'(\gt)|u|^2\exp(-\gm^2(\gt)t)\,dxdt}.
\end{equation}
If we suppose
\begin{equation}\label{2.9}
1-\frac{2\mu'(\gt)}{\mu^2(\gt)}\geq2^{-1},
\end{equation}
we derive from \eqref{2.8}
$$
f_{\mu,r}(\gt)+J_{\mu,r}(\gt)\leq-2\mu(\gt)^{-1}\frac{dJ_{\mu,r}(\gt)}{d\gt}.
$$
It is easy to check that this last inequality is equivalent to
$$
\frac{\mu(\gt)}2\exp\paran{\int_{\gt_1}^{\gt}\frac{\mu(s)}2\,ds}f_{\mu,r}(\gt)\leq
-\frac{d}{d\gt}\paran{J_{\mu,r}(\gt)\exp\paran{\int_{\gt_1}^\gt\frac{\mu(s)}2\,ds}}
\quad \forall\,\gt>\gt_1>k^{-1}.
$$
By integrating this inequality and using monotonicity of the function
$f_{\mu,r}(\gt)$ we get
$$
f_{\mu,r}(\gt_2)\int_{\gt_1}^{\gt_2}\frac{\mu(\gt)}2\exp\paran{\int_{\gt_1}^{\gt}\frac{\mu(s)}2\,ds}
d\gt+J_{\mu,r}(\gt_2)\exp\paran{\int_{\gt_1}^{\gt_2}\frac{\mu(s)}2\,ds}\leq
J_{\mu,r}(\gt_1)\quad \forall\, \gt_2>\gt_1>k^{-1}.
$$
Since
$$\frac{\mu(\gt)}2\exp\paran{\int_{\gt_1}^{\gt_2}
\frac{\mu(s)}2\,ds}=\frac{d}{d\gt}\paran{\exp\paran{\int_{\gt_1}^\gt\frac{\mu(s)}2\,ds}},$$
it follows from last the relation
\begin{equation}\label{2.10}
f_{\mu,r}(\gt_2)\sqbr{\exp\paran{\int_{\gt_1}^{\gt_2}\frac{\mu(s)}2\,ds}-1}+J_{\mu,r}(\gt_2)\exp
\paran{\int_{\gt_1}^{\gt_2}\frac{\mu(s)}2\,ds}\leq
J_{\mu,r}(\gt_1)\quad \forall\, \gt_2>\gt_1>k^{-1}.
\end{equation}
Now we have to define $\mu(\gt)$. Let $\vge>0$ and
\begin{equation}\label{2.11}
\mu(\gt)=\vge r^{-1}(\gt-k^{-1})\q \forall\,\gt>k^{-1}.
\end{equation}
One can easily verify that condition \eqref{2.9} is equivalent
to
\begin{equation}\label{2.12}
\gt\geq k^{-1}+2\vge^{-1/2}r^{1/2}.
\end{equation}
Now from \eqref{2.10} follow  two inequalities
\begin{multline}\label{2.13}
A(\gt_2):=\int_{Q^r(\gt_2)}\paran{|\nabla_xu|^2+\frac{\vge^2(\gt_2-k^{-1})^2}{r^2}|u|^2}dxdt
\leq
A(\gt_1)\\\times\exp\sqbr{-\frac{\vge\paran{(\gt_2-k^{-1})^2-(\gt_1-k^{-1})^2}}{4r}+\frac{\vge^2(\gt_2-k^{-1})}r}
\\\forall\, \gt_2>\gt_1>k^{-1}+2\vge^{-1/2}r^{1/2},
\end{multline}
and
\begin{multline}\label{New1}
f(r,\gt_2)\leq
A(\gt_1)\sqbr{\exp\paran{\frac{\vge\paran{(\gt_2-k^{-1})^2-(\gt_1-k^{-1})^2}}{4r}}-1}^{-1}
\exp\paran{\frac{\vge^2(\gt_2-k^{-1})^2}r} \\* \forall\,
\gt_2>\gt_1>k^{-1}+2\vge^{-1/2}r^{1/2}.
\end{multline}
In particular, for $\vge=8^{-1}$ we obtain from \eqref{2.13} and
\eqref{New1},
\begin{multline}\label{2.14}
\int_{Q^r(\gt)}\paran{|\nabla_xu|^2+\frac{(\gt-k^{-1})^2}{64r^2}|u|^2}dxdt\leq
e\exp\paran{-\frac{(\gt-k^{-1})^2}{64r}}\int_{Q^r(\gt_0^{(k)})}\paran{|\nabla_xu|^2+
\frac{|u|^2}{2r}}dxdt \\ \forall\,
\gt\geq\gt_0^{(k)}(r):=k^{-1}+4\sqrt{2}\sqrt{r},
\end{multline}
and
\begin{equation}\label{New2}
f(r,\gt)\leq\frac{e^2}{e-1}\exp\paran{-\frac{(\gt-k^{-1})^2}{64r}}\int_{Q^r(\gt_0^{(k)})}
\paran{|\nabla_xu|^2+\frac{u^2}{2r}}dxdt \quad \forall\,
\gt\geq\widetilde{\gt}_0^{(k)}(r):=k^{-1}+8\sqrt{r}.
\end{equation}
In order to have an estimate from above of the last factor in the
right-hand side of \eqref{2.14}, \eqref{New2}, we return to the
equation satisfied by $u$, multiply it by the test function
$u_k(x,t)\exp\paran{- t}$ and integrate over the domain
$Q^r=\BBR^N\times(0,r)$. As result of standard computations we
obtain, using \eqref{2.3},
\begin{multline}\label{2.15}
\int_{\BBR^N}|u_k(x,r)|^2\,dx+\int_{Q^r}\paran{|\nabla_xu_k|^2
+|u_k|^2+h(t)|u_k|^{q+1}}dxdt\\*\leq
\overline{c}\,\|u_{0,k}\|_{L_2(\BBR^N)}^2\leq cM_k\to\infty \text{
as }k\to\infty,\ \forall\,r\leq T.
\end{multline}
Due to \eqref{New2}, \eqref{2.15} it follows from \eqref{2.6}
\begin{multline}\label{2.16}
\int_{\BBR^N}|u(x,T)|^2\,dx+I_1(r)+I_2(r)+I_3(r)\\\leq
c_1\gt^{\frac{N(q-1)}{q+1}}h(r)^{-\frac{2}{q+1}}(-I'_2(r))^{\frac2{q+1}}
+c_2M_kr^{-1}\exp\paran{-\frac{(\gt-k^{-1})^2}{64r}}\quad
\forall\, \gt\geq\widetilde{\gt}_0^{(k)}(r).
\end{multline}
Relationships \eqref{2.14}, \eqref{New2} due to \eqref{2.15}
yield:
\begin{equation}\label{2.17}
f(r,\gt)+E_1(r,\gt)+\frac{(\gt-k^{-1})^2}{64r^2}E_2(r,\gt)\leq
c_2\,M_kr^{-1}\exp\paran{-\frac{(\gt-k^{-1})^2}{64r}}\quad
\forall\,\gt>\widetilde{\gt}_0^{(k)}(r).
\end{equation}
{\it Step 2. The first round of computations. } Next  we construct
some sequences $\{\gt_j\},\ \{r_j\},\ j=k,k-1,\ldots,1$. First we
explicit the choice of $M_k$ from condition \eqref{2.1}, let
namely
\begin{equation}\label{2.18}
M_k=e^{{e^k}}.
\end{equation}
Then we choose $\gt_k,\ r_k$ such that the following
relation is true,
\begin{equation}\label{2.19}
c_2\,r_k^{-1}\exp\paran{-\frac{\gt_k^2}{64r_k}}M_k=M_k^{\vge_0},\quad
0<\vge_0<e^{-1}
\end{equation}
where $c_2$ is from \eqref{2.16}, \eqref{2.17}. As consequence of
\eqref{2.19} and \eqref{2.18} we get
\begin{equation}\label{2.20}
\gt_k=8r_k^{1/2}\sqbr{(1-\vge_0)e^{k}+\ln r_k^{-1}+\ln
c_2}^{1/2}.
\end{equation}
In inequality \eqref{2.16} we fix $\gt=\gt_k+k^{-1}$, then due to
definition \eqref{2.19} it follows from \eqref{2.16},
\begin{multline}\label{2.21}
\int_{\BBR^N}|u(x,T)|^2\,dx+I_1(r)+I_2(r)+I_3(r)\\\leq
c_1(k^{-1}+\gt_k)^{\frac{N(q-1)}{q+1}}h(r)^{-\frac{2}{q+1}}(-I'_2(r))^{\frac2{q+1}}
+M_k^{\vge_0}\quad \forall\, r:0<r\leq r_k.
\end{multline}
$I_1(r),\ I_2(r),\ I_3(r)$ are nonincreasing functions which
satisfy, due to global a' priori estimate \eqref{2.15},
\begin{equation}\label{2.22}
I_1(0)+I_2(0)+I_3(0)\leq cM_k.
\end{equation}
Let us define the number $r_k$ by
\begin{equation}\label{2.23}
r_k=\sup\curlybr{r:I_1(r)+I_2(r)+I_3(r)\geq 2M_k^{\vge_0}}.
\end{equation}
Then it follows from  \eqref{2.21} the following differential
inequality
\begin{equation}\label{2.23*}
I_1(r)+I_2(r)+I_3(r)+\int_{\BBR^N}|u(x,T)|^2\,dx\leq 2c_1
(\gt_k+k^{-1})^{\frac{N(q-1)}{q+1}}h(r)^{-\frac{2}{q+1}}(-I'_2(r))^{\frac2{q+1}}\quad\forall\,
r\leq r_k.
\end{equation}
Solving it, we get
\begin{equation}\label{2.24}
I_1(r)+I_2(r)+I_3(r)\leq
c_3(\gt_k+k^{-1})^NH(r)^{-\frac{2}{q-1}}\quad \forall\,r\leq r_k,
\end{equation}
where
$$H(r)=\myint{0}{r}h(s)\,ds\quad\mbox{and }
c_3=\paran{\dfrac{2}{q-1}}^{2/(q-1)}\paran{2c_1}^{(q+1)/(q-1)}
$$
Next we will use more specific functions
$$
h(t)=\exp\paran{-\frac{\gw(t)}{t}},
$$
where $\gw(t)$ is  nondecreasing and satisfies the
following technical assumption
\begin{equation}\label{2.25}
t^{\ga_0}\leq\gw(t)\leq\gw_0=\text{const}\quad\forall\,
t:0<t<t_0,\ 0\leq\ga_0<1.
\end{equation}
It is easy to show by integration by parts the following
relation
$$
\int_0^r\exp\paran{-\frac{a\gw(t)}{t}}dt\geq\frac{1-\gd(r)}{(1-\ga_0)a}\cdot
\frac{r^2}{\gw(r)}\exp\paran{-\frac{a\gw(r)}{r}}\quad \forall\,
r>0,
$$
where $\gd(r)\to 0$ if $r\to0$. Therefore
\begin{equation}\label{2.26}
H(r)\geq\ovl{c}\frac{r^2}{\gw(r)}h(r),\ \ovl{c}=\text{const}>0.
\end{equation}
As a consequence we derive from \eqref{2.24}, using \eqref{2.20},
\begin{multline}\label{2.27}
I_1(r)+I_2(r)+I_3(r)\leq c_4\sqbr{8r_k^{\frac
12}\paran{(1-\vge_0)e^{k}+\ln r_k^{-1}+\ln c_2}^{\frac
12}+k^{-1}}^N\\\times\frac{\gw(r)^{\frac2{q-1}}}{r^{\frac4{q-1}}}\exp\paran{\frac{2\gw(r)}{(q-1)r}}\quad
\forall\,r\leq r_k.
\end{multline}
Comparing \eqref{2.23} and estimate \eqref{2.27} we deduce that
$r_k$ satisfies
\begin{equation}\label{2.28}
r_k\leq b_k,
\end{equation}
where $b_k$ is solution of equation
\begin{multline*}
c_4\sqbr{8b_k^{\frac 12}\paran{(1-\vge_0)e^{k}+\ln b_k^{-1}+\ln
c_2}^{\frac
12}+k^{-1}}^N\gw(b_k)^{\frac2{q-1}}b_k^{-\frac4{q-1}}\exp\paran{\frac{2\gw(b_k)}{(q-1)b_k}}\\=
2M_k^{\vge_0}=2\exp(\vge_0e^{k}).
\end{multline*}
This equation  may be rewritten in the form
\begin{multline}\label{2.29}
\ln c_4+\frac2{q-1}\ln
\paran{\frac{\gw(b_k)}{b_k}}+\frac2{q-1}\cdot\frac{\gw(b_k)}{b_k}\\
+N\ln\sqbr{8b_k^{\frac{N(q-1)-4}{2(q-1)N}}\paran{(1-\vge_0)\exp
k+\ln b_k^{-1}+\ln c_2}^{\frac
12}+k^{-1}b_k^{-\frac2{(q-1)N}}}=\ln2+\vge_0e^{k} \quad \forall\,
k\in\BBN.
\end{multline}
Since $s^{-1}\ln s\to 0$ as $s\to\infty$, it follows from equality
\eqref{2.29} that
\begin{multline}\label{2.30}
(1+c\gga(k))\vge_0e^{k}\geq
A_k+\frac2{q-1}\frac{\gw(b_k)}{b_k}\\:=N\ln\sqbr{8b_k^{\frac{N(q-1)-4}{2(q-1)N}}
\paran{(1-\vge_0)e^{k}+\ln b_k^{-1}+\ln c_2}^{\frac
12}+k^{-1}b_k^{-\frac2{N(q-1)}}}\\
+\frac2{q-1}\frac{\gw(b_k)}{b_k} \geq (1-\gga(k))\vge_0e^{k}\quad
\forall\, k\in\BBN,
\end{multline}
where $0<\gga(k)<1,\ \gga(k)\to0$ as $k\to\infty$. Keeping in
mind condition \eqref{2.25}, we obtain easily
\begin{equation}\label{2.31}
\frac{\gw (b_k)}{b_k}\geq b_k^{-(1-\ga_0)},\quad |A_k|\leq
c\paran{|\ln b_k|+k}\quad \forall\, k\in\BBN.
\end{equation}
Due to properties \eqref{2.31}, it follows from \eqref{2.30}
\begin{equation}\label{2.32}
ce^{k}>\frac{\gw (b_k)}{b_k}\geq d_1e^{k} \quad
\forall\,k\in\BBN,\ d_1>0.
\end{equation}
As a consequence of \eqref{2.32}, \eqref{2.31}  we obtain also
\begin{equation}\label{2.32*}
\ln b_k^{-1}\leq ck\quad \forall\, k\in\BBN.
\end{equation}
 Now using estimate \eqref{2.32} we are able to obtain
suitable upper estimate of $\gt_k$.  Thanks to
\eqref{2.28}, \eqref{2.32} and \eqref{2.32*} we deduce from \eqref{2.20}
$$
\gt_k\leq cb_k^{1/2}\exp\paran{\frac k2}\leq c\exp\paran{\frac k2}
\paran{\frac{\gw(b_k)}{d_1\exp
k}}^{1/2}=\frac{c}{d_1^{1/2}}\gw(b_k)^{1/2}.
$$
Using again estimate \eqref{2.32} and the monotonicity of the function
$\gw(s)$, we deduce from the above relation
\begin{equation}\label{2.33}
\gt_k\leq c\sqbr{\gw\paran{\frac{\gw_0}{d_1e^{k}}}}^{1/2}, \quad
\gw_0 \textrm{ is from (\ref{2.25})}.
\end{equation}
Therefore, from inequalities \eqref{2.17} and \eqref{2.27},
definitions \eqref{2.19}, \eqref{2.23} and property \eqref{2.28},
we derive the following estimates
\begin{equation}\label{2.34}
I_1(r_k)+I_2(r_k)+I_3(r_k)\leq 2M_k^{\vge_0}\quad \text{where
$r_k$ is from (\ref{2.28}), (\ref{2.23})},
\end{equation}
\begin{equation}\label{2.35}
f(r_k,\gt_k+k^{-1})+E_1(r_k,\gt_k+k^{-1})+\frac{\gt_k^2}{64r^2_k}E_2(r_k,\gt_k+k^{-1})
\leq M_k^{\vge_0},
\end{equation}
where $\gt_k$ is from (\ref{2.20}), (\ref{2.33}). Because
$\vge_0<e^{-1}$, it follows from definition \eqref{2.18} of
sequence $M_k$ that
\begin{equation}\label{2.36}
3M_k^{\vge_0}<cM_{k-1}\quad \forall\, k\geq k_0(c),
\end{equation}
where $c>0$ is arbitrary constant. Therefore, adding estimates
\eqref{2.34} and \eqref{2.35}, we obtain thanks to \eqref{2.36}
and the fact that $\gt_k\gg r_k$ (which follows from \eqref{2.19}), the
inequality
\begin{equation}\label{2.37}
f(r_k,\gt_k+k^{-1})+\sum\limits_{i=1}^3
I_i(r_k)+\sum\limits_{i=1}^2E_i(r_k,\gt_k+k^{-1})<cM_{k-1}\quad
\forall\, k\geq k_0(c).
\end{equation}
%
{\it Step 3. The second round of computations.} Next we introduce
the terms $r_{k-1},\ \gt_{k-1}$. Firstly we come back to
inequality \eqref{2.10}. Fixing here the  function
\begin{equation}\label{2.38}
\mu(t)=\vge r^{-1}(\gt-k^{-1}-\gt_k)\quad \forall\,
\gt>k^{-1}+\gt_k
\end{equation}
 instead of \eqref{2.11} and using estimates \eqref{2.12}--\eqref{New2}, we
obtain
\begin{equation}\label{2.39}\BA {l}
\displaystyle \int_{Q^r(\gt)}\paran{|\nabla_xu|^2+
\frac{(\gt-k^{-1}-\gt_k)^2|u|^2}{64r^2}}dxdt\\\phantom{---------}
\displaystyle\leq
e\exp\paran{-\frac{(\gt-k^{-1}-\gt_k)^2}{64r}}\int_{Q^r(\gt_0^{(k-1)}(r))}
\paran{|\nabla_xu|^2+
\frac{|u|^2}{2r}}dxdt\\[3mm]\phantom{--------------------}\forall\,
\gt>\gt_0^{(k-1)}(r):=k^{-1}+\gt_k+4\sqrt{2}\sqrt{r},
\EA\end{equation}
and
\begin{multline}\label{New3}
f(r,\gt)\leq\myfrac{e^2}{e-1}\exp\paran{-\myfrac{(\gt-k^{-1}-\gt_k)^2}{64r}}
\myint{Q^r(\gt_0^{(k-1)}(r))}{}
\paran{|\nabla_xu|^2+
\myfrac{|u|^2}{2r}}dxdt\\
 \forall\,
\gt\geq\widetilde{\gt}_0^{(k-1)}:=k^{-1}+\tau_k+8\sqrt{r}.
\end{multline}
The integral term in the right-hand side of \eqref{2.39}, \eqref{New3} is
estimated now by using estimate \eqref{2.37} obtained in the first
round of computation. So, we have
\begin{multline}\label{2.40}
\int_{Q^r(\gt_0^{(k-1)}(r))}
\paran{|\nabla_xu|^2+
\frac{u^2}{2r}}dxdt\leq (2r)^{-1}\sqbr{\sum\limits_{i=1}^3
I_i(r_k)+\sum\limits_{i=1}^2E_i(r_k,\gt_k+k^{-1})}\leq
c(2r)^{-1}M_{k-1}\\ \forall\, k>k_0(c),\ \forall\,r\geq r_k.
\end{multline}
Using this estimate we deduce from \eqref{2.39} and \eqref{New3}
\begin{multline}\label{2.41}
f(r,\gt)+E_1(r,\gt)+\frac{(\gt-\gt_k-k^{-1})^2}{64r^2}E_2(r,\gt)
\leq
c_2r^{-1}M_{k-1}\exp\paran{-\frac{(\gt-\gt_k-k^{-1})^2}{64r}}\\
\forall\,\gt\geq\widetilde{\gt}_0^{(k-1)}(r).
\end{multline}
This estimate is similar to estimate \eqref{2.17} from first
round. Now we have to deduce the analogue of estimate \eqref{2.24}. For
this we return to the starting relation \eqref{2.6}, where we now
estimate last term in right-hand side by estimate \eqref{New3},
using additionally \eqref{2.40}. As a result we have
\begin{multline}\label{2.46}
\sum\limits_{i=1}^3I_i(r)\leq
c_1\gt^{\frac{N(q-1)}{q+1}}h(r)^{-\frac2{q+1}}(-I'_2(r))^{\frac2{q+1}}+
c_2M_{k-1}r^{-1}\exp\paran{-\frac{(\gt-\gt_k-k^{-1})^2}{64r}}\\
\forall\,r\geq r_k, \ \forall\,
\gt\geq\widetilde{\gt}_0^{(k-1)}(r),
\end{multline}
which is analogous of estimate \eqref{2.16} from first round. Next we define
 the numbers $\gt_{k-1}$ and $r_{k-1}$ by  inequalities analogous to \eqref{2.20} and
\eqref{2.23},
\begin{equation}\label{2.47}
c_2r_{k-1}^{-1}M_{k-1}\exp\paran{-\frac{\gt_{k-1}^2}{64r_{k-1}}}=M_{k-1}^{\vge_0},\quad
0<\vge_0<e^{-1}
\end{equation}
\begin{equation}\label{2.48}
r_{k-1}=\sup\{r:I_1(r)+I_2(r)+I_3(r)\geq 2M_{k-1}^{\vge_0}\}.
\end{equation}
Now combining inequalities \eqref{2.23*} and \eqref{2.36}, and
using definitions \eqref{2.47}, \eqref{2.48}, we obtain the
following differential inequality
\begin{equation}\label{2.49} \sum\limits_{i=1}^{3}I_i(r)\leq
2c_1(\gt_{k-1}+\gt_k+k^{-1})^{\frac{N(q-1)}{q+1}}h(r)^{-\frac2{q+1}}(-I'_2(r))^{\frac2{q+1}}
\quad \forall\, r\leq r_{k-1}.
\end{equation}
Solving this differential inequality, we obtain an estimate similar to
\eqref{2.24}. Using property \eqref{2.26} we arrive to
\begin{equation}\label{2.50}
\sum\limits_{i=1}^{3}I_i(r)\leq c_4(\gt_{k-1}+\gt_k+k^{-1})^N\,
\frac{\gw(r)^{\frac2{q-1}}}{r^{\frac4{q-1}}}
\exp\paran{\frac{2\gw(r)}{(q-1)r}}\quad \forall\, r\leq r_{k-1}.
\end{equation}
As in first round we express from \eqref{2.47} $\gt_{k-1}$ as
function $\gt_{k-1}(r_{k-1})$ (the analogue of \eqref{2.20})
\begin{equation}\label{2.51}
\gt_{k-1}=8r_{k-1}^{1/2}[(1-\vge_0)\exp(k-1)+\ln r_{k-1}^{-1}+\ln
c_2]^{1/2}.
\end{equation}
Inserting this expression of $\gt_{k-1}$ into \eqref{2.50} and
then comparing the obtained inequality with definition \eqref{2.48}, we
deduce an estimate similar to \eqref{2.28},
\begin{equation}\label{2.52}
r_{k-1}\leq b_{k-1},
\end{equation}
where $b_{k-1}$ is solution of equation
\begin{multline}\label{2.53}
c_4\sqbr{8b_{k-1}^{1/2}\paran{(1-\vge_0)\exp(k-1)+\ln b_k^{-1}+\ln
c_2}^{1/2}+\gt_k+k^{-1}}^N\\\times\frac{\gw(b_{k-1})^{\frac2{q-1}}}{b_{k-1}^{\frac4{q-1}}}
\exp\paran{\frac{2\gw(b_{k-1})}{(q-1)b_{k-1}}}=2M_{k-1}^{\vge_0}=2\exp(\vge_0\exp(k-1)).
\end{multline}
From \eqref{2.41}, and  due to definition \eqref{2.47}, it follows
\begin{equation}\label{2.54}
f(r_{k-1},\gt_{k-1}+\gt_k+k^{-1})+\frac{\gt_{k-1}^2}{64r_{k-1}}E_2(r_{k-1},\gt_{k-1}+\gt_k+k^{-1})
+E_1(r_{k-1},\gt_{k-1}+\gt_k+k^{-1})\leq M_{k-1}^{\vge_0}.
\end{equation}
From \eqref{2.50}, due to \eqref{2.51}, \eqref{2.52}, \eqref{2.53},
it follows
\begin{equation}\label{2.55}
I_1(r_{k-1})+I_2(r_{k-1})+I_3(r_{k-1})\leq2M_{k-1}^{\vge_0}.
\end{equation}
Summing \eqref{2.54}, \eqref{2.55} and using property
\eqref{2.36}, we deduce new global {\it a priori} estimate (the analogous of
\eqref{2.37}) which is the main starting information for the next round
of computation
\begin{equation}\label{2.56}
f(r_{k-1},\gt_{k-1}+\gt_k+k^{-1})+\sum\limits_{i=1}^3I_i(r_{k-1})+\sum\limits_{i=1}^2
E_i(r_{k-1},\gt_{k-1}+\gt_k+k^{-1})\leq cM_{k-2}.
\end{equation}
We are ready now for the next round of computations, introducing
the function
$$
\mu(t)=\vge r^{-1}(\gt-k^{-1}-\gt_k-\gt_{k-1}) \quad
\forall\,\gt>k^{-1}+\gt_k+\gt_{k-1}
$$
instead of \eqref{2.38} and estimate \eqref{2.56} instead of
\eqref{2.37}. We realize $j$ rounds of such computations. As
result we obtain
\begin{equation}\label{2.57}
f\paran{r_{k-j},
\sum_{l=0}^j\gt_{k-l}+k^{-1}}+\sum_{i=1}^3I_i(r_{k-j})+\sum_{i=1}^2
E_i\paran{r_{k-j}, \sum_{l=0}^j\gt_{k-l}+k^{-1}}\leq cM_{k-j-1},
\end{equation}
which was our main aim.

{\it Step 4. The control of $r_{k-j}, \sum_{l=0}^j\gt_{k-l}$ as
$j\to k$ with arbitrary $k\in\BBN$}. It is clear that $r_{k-j}, \
\gt_{k-j}$ are defined by the conditions (see \eqref{2.47},
\eqref{2.48})
\begin{equation}\label{2.58}
c_2r_{k-j}^{-1}M_{k-j}\exp\paran{-\frac{\gt_{k-j}^2}{64r_{k-j}}}=M_{k-j}^{\vge_0},
\quad 0<\vge_0<e^{-1}.
\end{equation}
\begin{equation}\label{2.59}
r_{k-j}=\sup\bigl\{r:I_1(r)+I_2(r)+I_3(r)\geq
2M_{k-j}^{\vge_0}\bigr\}.
\end{equation}
Similarly to \eqref{2.51}--\eqref{2.53} we deduce that
\begin{equation}\label{2.60}
\gt_{k-j}=8r_{k-j}^{1/2}\sqbr{(1-\vge_0)e^{k-j}+\ln
r_{k-j}^{-1}+\ln c_2}^{1/2},
\end{equation}
\begin{equation}\label{2.61}
r_{k-j}\leq b_{k-j},
\end{equation}
$ \text{ where } b_{k-j} \text{ satisfies  }$
\begin{multline}\label{2.62}
c_4\sqbr{8b_{k-j}^{1/2}\paran{(1-\vge_0)e^{k-j}+\ln
b_{k-j}^{-1}+\ln
c_2}^{1/2}+\sum_{i=0}^{j-1}\gt_{k-i}+k^{-1}}^N\\\times
\frac{\gw(b_{k-j})^{\frac2{q-1}}}{b_{k-j}^{\frac4{q-1}}}
\exp\paran{\frac{2\gw(b_{k-j})}{(q-1)b_{k-j}}}=2M_{k-j}^{\vge_0}=2\exp(\vge_0e^{k-j}).
\end{multline}
In the first round of computations we have obtained the upper estimate
\eqref{2.33} for $\gt_k$. Let us suppose by induction that the
following estimate is true
\begin{equation}\label{2.63}
\gt_{k-i}\leq
c\sqbr{\gw\paran{\frac{\gw_0}{d_1\exp(k-i)}}}^{1/2}\quad
\forall\,i\leq j-1.
\end{equation}
We have to prove that estimate \eqref{2.63} holds also for $i=j$. Obviously condition \eqref{2.62} is equivalent to
(see \eqref{2.29})
\begin{equation}\label{2.64}
\ln
c_4+\frac{2}{q-1}\ln\paran{\frac{\gw(b_{k-j})}{b_{k-j}}}+\frac2{q-1}
\cdot \frac{\gw(b_{k-j})}{b_{k-j}}+A_k^{(j)}=\ln2+\vge_0e^{k-j},
\end{equation}
where
$$A_k^{(j)}=
N\ln\sqbr{b_{k-j}^{\frac{N(q-1)-4}{2(q-1)N}}
\paran{(1-\vge_0)e^{k-j}+\ln(b_{k-j}^{-1})+\ln c_2}^{1/2}+
\frac{k^{-1}+\sum\limits_{i=0}^{j-1}\gt_{k-i}}{b_{k-j}^{\frac2{(q-1)N}}}}.
$$
Because of the induction assumption \eqref{2.63}
$$
\sum_{i=0}^{j-1}\gt_{k-i}\leq
c\sum_{i=0}^{j-1}\sqbr{\gw\paran{\frac{\gw_0}{d_1\exp(k-i)}}}^{1/2}\leq
c\int_0^1\frac{\gw(s)^{1/2}}s\, ds: =cL,
$$
therefore
\begin{equation}\label{2.65}
|A_k^{(j)}|\leq c \paran{|\ln b_{k-j}|+(k-j)+\ln L}.
\end{equation}
From \eqref{2.64} due to  \eqref{2.65} we derive easily
\begin{equation}\label{2.66}
ce^{k-j}\geq\frac{\gw(b_{k-j})}{b_{k-j}}\geq d_1e^{k-j}\quad
\forall\,j: k-j\geq k_0=k_0(L),
\end{equation}
where $k_0<\infty$ do not depend on $k$. From \eqref{2.66} it
follows in particular
\begin{equation}\label{2.67}
\ln b_{k-j}^{-1}\leq c(k-j)\quad \forall\,j: k-j\geq k_0.
\end{equation}
Thanks to \eqref{2.61} and
properties \eqref{2.66}, \eqref{2.67},  we derive from \eqref{2.60},
\begin{multline}\label{2.68}
\gt_{k-j}\leq 8b_{k-j}^{1/2}\paran{(1-\vge_0)e^{k-j}+\ln
b_{k-j}^{-1}+\ln c_2}^{1/2}\\\leq
cb_{k-j}^{1/2}\exp\paran{\frac{k-j}2}\leq \frac
c{d_1^{1/2}}[\gw(b_{k-j})]^{1/2}\quad \forall\,j: k-j\geq k_0(L).
\end{multline}
Using again estimate \eqref{2.66} and monotonicity of $\gw(s)$ we
deduce from \eqref{2.68}
\begin{equation}\label{2.69}
\gt_{k-j}\leq
c\sqbr{\gw\paran{\frac{\gw_0}{d_1e^{k-j}}}}^{1/2}\quad
\forall\,j: k-j\geq k_0(L).
\end{equation}
Thus, we have proved  by induction estimate \eqref{2.63}, for arbitrary
$k-j\geq k_0(L)$ with $r_i,\ \gt_i$ satisfying \eqref{2.61},
\eqref{2.62} and \eqref{2.69}. \medskip

\noindent {\it Step 5. Completion of the proof. }We fix now $n>k_0(L)$ and take $j=k-n$
in \eqref{2.57}. This leads to
\begin{equation}\label{2.70}
f\paran{r_n,\sum_{l=0}^{k-n}\gt_{k-l}+k^{-1}}+\sum_{i=1}^3I_i(r_n)+\sum_{i=1}^2
E_i\paran{r_n,\sum_{l=0}^{k-n}\gt_{k-l}+k^{-1}}\leq cM_{n-1}\quad
\forall\, n>k_0(L).
\end{equation}
Next we have
\begin{equation}\label{2.71}
\sum_{l=0}^{k-n}\gt_{k-l}\leq \sum_{i=n}^{\infty}\gt_i\leq
c\sum_{i=n}^{\infty}\sqbr{\gw\paran{\frac{\gw_0}{d_1\exp
i}}}^{1/2}\leq
c\int_0^{\frac{\gw_0}{d_1\exp(n-1)}}\frac{\gw(s)^{1/2}}s\,ds\to0\text{
as }n\to\infty.
\end{equation}
Therefore, for arbitrary small $\gd>0$, we can find  and fix
$n=n(\gd)<\infty$ such that from \eqref{2.70} follows uniform with
respect to $k\in\BBN$ {\it a priori} estimate,
\begin{equation}\label{2.72}
\sup_{t>0}\int_{|x|>\gd}|u_k(x,t)|^2\,dx
+\int_0^T\int_{|x|>\gd}\paran{|\nabla_xu_k|^2+|u_k|^2}dxdt\leq
C=C(\gd)<\infty \quad \forall\, k\in\BBN.
\end{equation}

Since $u_k(x,0)=0\ \ \forall\,|x|>k^{-1}\ \ \forall\,k\in\BBN$, it
follows from  \eqref{2.72} that $u_\infty(x,0)=0\ \forall\,x\ne0$, which ends the proof.
\qeda \\[3mm]

\section{Regional initial blow-up for equation with \\exponential absorption.}
\setcounter{equation}{0}

The local energy method we have used in the proof of Theorem
\ref{t:locth1} is based on the sharp interpolation theorems for
functional Sobolev spaces, which are natural tool for the study of
solutions of equations with power nonlinearities. Here we propose
the adaptation of mentroned method to the equations with nonpower
nonlinearities.

Thus, we consider the Cauchy problem
\begin{equation}
\label{hexp-u}
\left\{\BA {l}
\partial_{t}u-\Gd u+h(t)(e^{u}-1)=0\quad \mbox {in }Q^{\infty}\\[2mm]
u(x,0)=k\gd_{0},
\EA\right.
\end{equation}

\bth{locth2} Assume $h(t)=e^{-e^{\gw(t)/t}}$ where $\gw\in
C([0,\infty))$ satisfies the same asumptions as in \rth {locth1}.
Then solution $u_{k}$ always exists and
$u_{\infty}:=\lim_{k\to\infty}u_{k}$ has a point-wise singularity
at $(0,0)$. \es


\noindent \Proof We will consider the family $u_k(x,t)$ of
solutions of regularized problems:
\begin{equation}\label{3.0}\left\{\BA {l}
u_t-\Gd u+h(t)(e^{u}-1)=0\quad\mbox{in } Q^T,\\[2mm]
u(x,0)=u_{0,k}(x)=M^{1/2}_kk^{- N/2}\gd_k(x)\forevery x\in\BBR^N,
\EA\right.\end{equation} where $\gd_k$ is nonnegative, continuous
with compact support in $B_{k^{-1}}$, satisfies estimate
\eqref{2.3} and converges weakly to $\gd_{0}$ as $k\to\infty$,
$\curlybr{M_k}$ satisfies condition \eqref{dini1}. Let us
introduce the energy functions (we omit index $k$ in $u_k$):
\begin{equation}\label{3.1}
I_{1,0}(r)=\int_{Q_r}|\nabla_xu|^2\,dxdt,\quad
I_q(r)=(q!)^{-1}\int_{Q_r}h(t)|u|^{q+1}\,dxdt, \quad
I_{3,0}(r)=\int_{Q_r}|u|^2\,dxdt.
\end{equation}
Multiplying \eqref{3.0} by
$u(x,t)\exp\paran{-\dfrac{t-r}{1+T-r}}$, integrating in $Q_r$ and
using equality
$$s(e^s-1)=\sum\limits_{q=1}^\infty\dfrac{s^{q+1}}{q!},$$
we obtain
easily
\begin{multline}\label{3.2}
I_{1,0}(r)+\sum_{l=1}^{\infty}I_l(r)+I_{3,0}(r)\leq
c(q!)^{2/(q+1}\gt^{N(q-1)/(q+1)}h(r)^{-2/(q+1)}(-I'_q(r))^{2/(q+1)}
\\+ c\int_{\Gw(\gt)}|u(x,r)|^2\,dx \quad \forall\,\gt>0,\
\forall\,r:0<r<T,\ \forall\,q\in\BBN.
\end{multline}
We introduce the additional energy functions
\begin{equation}\label{3.3}
f(r,\gt)\text{ from (\ref{2.6*})},\quad
E_{1,0}(r,\gt)=\int_{Q^r(\gt)}|D_xu|^2\,dxdt,\quad
E_{2,0}(r,\gt)=\int_{Q^r(\gt)}|u|^2\,dxdt.
\end{equation}
Instead of \eqref{2.15} we derive the following global
{\it a priori} estimate:
\begin{multline}\label{New4}
\int_{\mathbb{R}^N}|u_k(x,r)|^2\,dx+\int_{Q^r}\paran{|\nabla_xu|^2+|u_k|^2+h(t)\sum_{l=1}^\infty
\frac{|u_k|^{l+1}}{l!}}dxdt\\\leq \overline{c}\,\| u_{0,k}
\|^2_{L_2(\mathbb{R}^N)}\leq cM_k \quad \forall\,r<T.
\end{multline}
Using estimate \eqref{New4} instead of \eqref{2.15} in a similar way
as in the proof of \rth {locth1}, we obtain the following  inequality, analogous to
\eqref{2.17},
\begin{multline}\label{3.4}
f(r,\gt)+E_{1,0}(r,\gt)+\frac{(\gt-k^{-1})^2}{64r^2}E_{2,0}(r,\gt)+\\\leq
c_2M_kr^{-1}\exp\paran{-\frac{(\gt-k^{-1})^2}{64r}}\quad \forall\,
\gt\geq\widetilde{\gt}_0^{(k)}(r)=k^{-1}+8\sqrt{r}.
\end{multline}
Using this estimate we deduce from \eqref{3.2}
\begin{multline}\label{3.5}
I_{1,0}(r)+\sum_{l=1}^\infty I_l(r)+I_{3,0}(r)\leq
c(q!)^{\frac2{q+1}}\gt^{\frac{N(q-1)}{q+1}}h(r)^{-\frac2{q+1}}(-I'_q(r))^{\frac2{q+1}}\\
+c_2M_kr^{-1}\exp\paran{-\frac{(\gt-k^{-1})^2}{64r}}\quad
\forall\,\gt\geq\widetilde{\gt}_0^{(k)}(r),\ \forall\, q\in\BBN.
\end{multline}
Next, we define the numbers $\gt_k,\ r_k$. Firstly, set
\begin{equation}\label{3.6}
r_k:=\sup\curlybr{r:I_{1,0}(r)+\sum_{l=1}^\infty
I_l(r)+I_{3,0}\geq 2M_k^{\vge_0}},\quad 0<\vge_0<e^{-1}.
\end{equation}
Then we fix the sequence $\{M_k\}$ by \eqref{2.18} again and $\gt_k$
by inequalities \eqref{2.19}, \eqref{2.20}. Thanks to these
definitions we derive the following
series of inequalities from relations \eqref{3.5}
\begin{equation}\label{3.7}
I_{1,0}(r)+\sum_{l=1}^\infty I_l(r)+I_{3,0}(r)\leq
2c_1(q!)^{\frac2{q+1}}(\gt_k+k^{-1})^{\frac{N(q-1)}{q+1}}h(r)^{-\frac2{q+1}}(-I'_q(r))^{\frac{2}{q+1}}
\quad \forall\, q\in\mathbb{N}, \ \forall\, r\leq r_k.
\end{equation}
Solving these differential inequalities we obtain the
estimates
\begin{equation}\label{3.8}
I_{1,0}(r)+\sum_{l=1}^\infty I_l(r)+I_{3,0}(r)\leq
c_3(\gt_k+k^{-1})^N(q!)^{\frac2{q-1}}H(r)^{-\frac2{q-1}}\quad
\forall\, r\leq r_k,\ \forall\,q\in\mathbb{N},
\end{equation}
where $H(r)$ is from \eqref{2.24}. We have now to optimize
estimate \eqref{3.8} with respect  to parameter $q$. By integration by parts,  it is easy to
check the following inequality
\begin{equation}\label{3.9}
H(r)\geq
\overline{c}\frac{r^2}{\gw(r)}\exp\paran{-\frac{\gw(r)}{r}}h(r)\quad
\forall\, r>0,\ \overline{c}>0.
\end{equation}
Using Stirling formula $q!\sim\paran{\dfrac qe}^q$ and estimate
\eqref{3.9}, we deduce from \eqref{3.8}
\begin{equation}\label{3.10}
I_{1,0}(r)+\sum_{l=1}^\infty I_l(r)+I_{3,0}(r)\leq
c_4(\gt+k^{-1})^NF_q(r) \quad \forall\, r\leq r_k,
\end{equation}
where
$$
F_q(r)=q^2\gw(r)^{\frac2{q-1}}r^{-\frac4{q-1}}\exp\paran{\frac2{q-1}\cdot\frac{\gw(r)}{r}}
\exp\sqbr{\frac2{q-1}\exp\paran{\frac{\gw(r)}r}}.
$$
Fixing here the optimal value of the parameter $q$:
$$
q=\widetilde{q}:=\sqbr{2\exp\paran{\frac{\gw(r)}r}},
$$
where $[a]$ denotes the enteger part of $a$, we obtain easily
$$
F_{\widetilde{q}}\leq c\exp\paran{\frac{2\gw(r)}r}.
$$
Therefore it follows from \eqref{3.10},
\begin{equation}\label{3.11}
I_{1,0}(r)+\sum_{l=1}^\infty I_l(r)+I_{3,0}(r)\leq
c_5(\gt_k+k^{-1})^N\exp\paran{\frac{2\gw(r)}r}\quad \forall\,r\leq
r_k.
\end{equation}
Comparing now definition \eqref{3.6} of $r_k$ and estimate
\eqref{3.11}, and  using additionally the expression
\eqref{2.20} of $\gt_k$, we obtain
\begin{equation}\label{3.12}
r_k\leq b_k,
\end{equation}
where $b_k$ is defined by the equation
\begin{multline}\label{3.13}
c_5\sqbr{8b_k^{1/2}((1-\vge_0)e^{k}+\ln b_k^{-1}+\ln
c_2)^{1/2}+k^{-1}}^N\exp\paran{\frac{2\gw(b_k)}{b_k}}\\=2M_k^{\vge_0}=2\exp(\vge_0\exp
k), \quad 0<\vge_0<e^{-1}.
\end{multline}
By an analysis similar to Step 2 in the proof of \rth {locth1}, we obtain estimates
\eqref{2.30}--\eqref{2.32*} for $b_k$. Then we prove the validity of
estimate \eqref{2.33} for $\gt_k$. As a consequence of estimates
\eqref{3.4}, \eqref{3.11},  thanks to to definitions \eqref{2.20},
\eqref{3.6} of $\gt_k,\ r_k$ and the previous estimates of $\gt_k,\
r_k$, we get
\begin{align*} &I_{1,0}(r)+\sum_{l=1}^\infty
I_l(r)+I_{3,0}(r)\leq2 M_k^{\vge_0},\\
&f(r_k,\gt_k+k^{-1})+E_{1,0}(r_k,\gt_k+k^{-1})+\frac{\gt_k^2}{64r_k^2}E_{2,0}(r_k,\gt_k+k^{-1})\leq
M_k^{\vge_0}.
\end{align*}
Summing these inequalities, and using definition of $\{M_k\}$  and
property $\gt_k\gg r_k$, we obtain an analogue of estimate \eqref{2.37}, namely,
\begin{equation}\label{3.14}
f(r_k,\gt_k+k^{-1})+I_{1,0}(r_k)+\sum_{l=1}^\infty
I_l(r_k)+I_{3,0}(r_k)+E_{1,0}(r_k,\gt_k+k^{-1})+E_{2,0}(r_k,\gt_k+k^{-1})\leq
cM_{k-1}.
\end{equation}
Using \eqref{3.14} as global {\it a priori} estimate instead of
\eqref{New4} and providing a second round of computations similar
to \eqref{2.38}--\eqref{2.52} we derive a second global {\it a priori}
estimate analogous to \eqref{2.56},
\begin{multline*}
f(r_{k-1},\gt_{k-1}+\gt_k+k^{-1})+I_{1,0}(r_{k-1})+\sum_{l=1}^\infty
I_l(r_{k-1})+I_{3,0}(r_{k-1})\\+E_{1,0}(r_{k-1},\gt_{k-1}+\gt_k+k^{-1})+E_{2,0}(r_{k-1},\gt_{k-1}+
\gt_k+k^{-1})\leq cM_{k-2}.
\end{multline*}
Repeating such rounds $j$-times we derive a corresponding  analogue of
relation \eqref{2.57}. It is easy to see that estimate
\eqref{2.71} for constructed shifts $\gt_{k-i}$ remains
valid. This fact, similar to what was used in the proof of \rth {locth1}, yields to the conclusion.\qeda

\section{The porous media equation with absorption}
\setcounter{equation}{0}

In this section we consider the following problem dealing with
fundamental solutions of the porous media equation with time
dependent absorption,
\begin{equation}
\label{pme-k}\left\{\BA{l} \prt_{t} u-\Gd( |u|^{m-1}u)
+h(t)|u|^{q-1}u=0\quad \mbox {in }Q^T\\[2mm]
u(x,0)=k\gd_0. \EA\right.
\end{equation}
It is standard to assume that $h\geq 0$ is a continuous function
and $m,q$ are positive real numbers. By a solution we mean a
function $u\in L_{loc}^1( Q^T)$ such that $u^m\in L_{loc}^1(
Q^T)$, $hu^q\in L_{loc}^1( Q^T)$ and
\begin{equation}
\label {weak*}
\dint_{Q^T}\left(-u\prt_t\gf-|u|^{m-1}u\Gd\gf+h(t)|u|^{q-1}u\gf\right)dxdt
=k\gf(0,0)
\end {equation}
for any $\gf\in C_0^{2,1}(\BBR^N\ti [0,T))$.
If $h\equiv 0$ and $m>(N-2)_+/N$ this problem admits a solution for any $k>0$. When $m>1$ this solution has the following form
\begin{equation}\label{barb}
B_k(x,t)=t^{-\ell}\left(C_k-\myfrac{(m-1)\ell}{2mN}\myfrac{\abs x^2}{t^{2\ell/N}}\right)^{1/(m-1)}_+,
\end {equation}
where
\begin{equation}\label{barb2}
\ell=\myfrac{N}{N(m-1)+2}\quad\mbox{and } C_k=a(m,N)k^{2(m-1)\ell/N}.
\end {equation}
Since  $B_k$ is a supersolution for problem (\ref{pme-k}), a sufficient condition for existence (and uniqueness) of $u_k$ is
\begin{equation}\label {Bar-q}
\dint_{Q^T}B^q_k(x,t)h(t)dxdt<\infty.
\end{equation}
By the change of variable $y=t^{\ell/N}x$ this condition is independent of $k>0$ and we have
\bprop{fundpm} Assume $m>1$, $q>0$. If
\begin{equation}\label {Bar-q1}
\int_0^1h(t)t^{\ell-\ell q}dt<\infty,
\end{equation}
then problem (\ref{pme-k}) admits a unique positive solution $u=u_k$. In the particular case where $h(t)=O(t^{\ga})$ ($\ga\geq 0$), the condition is
\begin{equation}\label {Bar-q2}
\ga>\myfrac{N(q-m)-2}{N(m-1)+2}.
\end{equation}
\es

We recall that if  $q>1$ and $m>(N-2)_+/N$, any solution of the porous media equation with absorption is bounded from above by the maximal solution $U_h$ expressed by
\begin{equation}\label {max1}
U_h(t)=\left((q-1)\myint {0}{t}h(s)\,ds\right)^{-1/(q-1)}.
\end{equation}
\bth{pmth1} Assume $q+1>2m>2$ and $h\in C((0,\infty))$ is nondecreasing, positive and satisfies
$h(t)=O(t^{(q-m)/(m-1)})$ as $t\to 0$. Then for any $k>0$ $u_k$ exists and $\lim_{k\to\infty}u_k:=u_{\infty}=U_h$.
\es
\Proof We first notice that
 $$q+1>2m>2\Longrightarrow q>m>1\quad\mbox {and } \;
\myfrac{q-m}{m-1}>\myfrac{N(q-m)-2}{N(m-1)+2}.$$
{\it Step 1. Case $q<m+2/N$. } In this range of value we know \cite {PT} that there exists a nonnegative very singular solution $v=v_\infty$ to
\begin{equation}\label{pme3}
\prt_tv-\Gd v^m+v^q=0\quad\mbox {in }\;Q^T,
\end{equation}
and $v_\infty=\lim_{k\to°}v_k$, where the $v_k$ are solutions of the same equation with initial data $k\gd_0$. Furthermore, $v_\infty$ is unique \cite {KV}, radial with respect to $x$ and has the following form
$$v_\infty(x,t)=t^{-1/(q-1)}F(\abs x/t^{(q-m)/2(q-1)}),
$$
where $F$ solves
\begin{equation}\label{pme4} \left\{\BA {l}
(F^m)''+\myfrac{N-1}{\eta}(F^m)'+\myfrac{q-m}{2(q-1)}\eta F'+\myfrac {1}{q-1}F-F^q=0\quad\mbox {in }\;(0,°)\\[3mm]
F'(0)=0\;\mbox {and }\;\lim_{\eta\to\infty}\eta^{2/(q-m)}F(\eta)=0.
\EA\right.\end{equation}
Actually $F$ has compact support in $[0,\xi_0]$ for some $\xi_0>0$. Let $\gga=(q-m)/(m-1)$, then for any $\ge >0$, $u=u_\infty$ satisfies, for some $c>0$,
$$\prt_tu-\Gd u^m+c\ge^{\gga}u^q\geq 0\quad\mbox {in }\;Q^\ge.
$$
If we set $w_\ge(x,t)=a^\gth v_\infty(x,at)$ with $\gth=1/(m-1)-$ and
$a=\ge^{-1}c^{-(q-1)/(q-m)}$, then
$$\prt_tw_\ge-\Gd w_\ge^m+c\ge^{\gga}w_\ge^q= 0\quad\mbox {in }\;Q^T.
$$
\medskip
By comparison $u_\infty\geq w_\ge$ in $Q^\ge$. If we take in particular $t=\ge$, it implies
\begin{equation}\label{pme5}
u_\infty(x,t)\geq c^{-1/(q-m)}
t^{-1/(m-1)}v_\infty(x,c^{-(m-1)/(q-m)}) =c^{-1}t^{-1/(m-1)}
F(c^{(m-1)/2(q-1)}\abs x) \end{equation} If $\abs
x<\xi_c=c^{-(m-1)/2(q-1)}\xi_0$, we derive that $\lim_{t\to
0}u_\infty(x,t)=\infty$, locally uniformly in $B_{\xi_c}$. This
implies $u_{\infty}=U_h$.\medskip

\noindent {\it Step 2. Case $q\geq m+2/N$. } We give an alternative proof valid for all $q$. We first observe that it is sufficient to prove the result when $h(t)$ is replaced by $t^{\gga}$. If we look for a family of transformations $u\mapsto T_\ell(u)$ under the form
$$ T_\ell(u)(x,t)=\ell^{\ga}u(\ell^\gb x,\ell t)\forevery (x,t)\in Q^\infty,\;\forall\ell>0
$$
which leaves the equation
\begin{equation}\label{pme6}\BA {l}
\prt _tu-\Gd |u|^{m-1}u+t^{\gga}|u|^{q-1}u=0
\EA\end{equation}
invariant, we find $\ga=(1+\gga)/(q-1)$ and $\gb=(q-m-\gga (m-1))/2(q-1)$. Due to the value of $\gga$, we have $\gb=0$. Because of uniqueness and the value of the initial mass
\begin{equation}\label{pme7}\BA {l}
T_\ell(u_k)=u_{\ell^\ga k}\forevery \ell>0,\;\forall k>0\Longrightarrow T_\ell(u_\infty)=u_{\infty}\forevery \ell>0.
\EA\end{equation}
Therefore
$$\ell^\ga u_\infty(x,\ell t)=u_\infty(x, t)\forevery (x,t)\in Q^\infty,\;\forall\ell>0.
$$
In particular, if we take $\ell=t^{-1}$,
$$u_\infty(x, t)=t^{-\ga} u_\infty(x,1)=t^{-\ga}\gf(x).
$$
Plugging this decomposition into (\ref{pme6}) yields to
$$-\ga t^{-\ga-1}\gf-t^{-\ga m}\Gd \gf^m+t^{\gga-\ga q}\gf^q=0,
$$
where all the exponents of $t$ coincide since
$$\ga m=\myfrac{m}{m-1}\;,\;\ga q-\gga =\myfrac{m}{m-1}\;\mbox{ and }\;
\ga+1=\myfrac{m}{m-1}.
$$
Therefore $\gf$ is a positive and radial (as the $u_k$ are) solution of
$$-\ga \gf-\Gd \gf^m+\gf^q=0\quad\mbox {in }\;\BBR^N.
$$
Setting $\psi=\gf^m$ yields to
\begin{equation}\label{pme7}\BA {l}
-\Gd \psi-\myfrac{1}{m-1}\psi^{1/m}+\psi^{q/m}=0\quad\mbox {in }\;\BBR^N.
\EA\end{equation}
Clearly $\psi=\psi_0=(m-1)^{-m/(q-1)}$ is a solution. By a standard variation of the Keller-Osserman estimate, any solution is bounded from above by $\psi_0$. Putting
$\tilde \psi (x)=A\psi (a)$, it is easy to find $A>0$ and $a>0$ such that
\begin{equation}\label{pme7}\BA {l}
-\Gd \tilde\psi-\tilde\psi^{1/m}+\tilde\psi^{q/m}=0\quad\mbox {in }\;\BBR^N,
\EA\end{equation}
with $0\leq \tilde\psi\leq 1$. Writting $\tilde\psi$ as a solution of an ODE, we derive
$$\tilde\psi(r)=\tilde\psi(0)+
\myint{0}{r}s^{1-n}\myint{0}{s}(\tilde\psi^{q/m}-\tilde\psi^{1/m})\gs^{n-1}ds
\forevery r>0.$$
If $\tilde\psi^{q/m} $ is not constant with value $1$, the right-hand side of the above inequality is decreasing with respect to $r$, and the only possible nonnegative limit is $0$, by La Salle principle. Thus
$$\tilde\psi''+\myfrac{N-1}{r}\tilde\psi'+\myfrac{1}{2}\tilde\psi^{1/m}\leq 0
$$
for $r\geq r_0$, large enough. If $N=2$, we set $\gt=\ln r$, $\Psi(\gt)=\tilde\psi(r)$ and get
$$\Psi''+\myfrac{1}{2}e^{2\gt}\Psi^{1/m}\leq 0
$$
for $\gt\geq \ln r_0$. The concavity of $\Psi$ yields a contradiction. If $N\geq 3$, we set $\gt =r^{N-2}/(N-2)$ and $\Psi(\gt)=r^{N-2}\tilde\psi(r)$. Then $\Psi$ satisfies
$$\Psi''+c_N\gt^{(4-N)/(N-2)-1/m}\Psi^{1/m}\leq 0.
$$
Again the concavity yields a contradiction. In any case we obtain that $\Psi=1$, or, equivalently $\psi=\psi_0$ and finally, $u_\infty= t^{-1/(m-1)}\psi_0^{1/m}$. \qeda\\

\bth {pmth2} Assume $q>m>1$ and $h\in C((0,\infty))$ is nondecreasing, positive. If
$h(t)=t^{(q-m)/(m-1)}\gw^{-1}(t)$ with
$\gw(t)\to 0$ as $t\to 0$, and
\begin{equation}
\label{intr6*} \int_0^1\gw^\gth(s)\myfrac {ds}{s}<\infty,
\end{equation}
where
$$
\gth=\myfrac{m^2-1}{[N(m-1)+2(m+1)](q-1)},
$$
 then $u_{\infty}:=\lim_{k\to\infty}u_{k}$ has a point-wise singularity at $(0,0)$\es
\Proof The structure of the proof is similar to the one of \rth
{locth1}. We  study the asymptotic behaviour as $k\to\infty$ of
solutions $u=u_k(x,t)$ of the regularized Cauchy problem
\begin{equation}\label{4.1}\left\{\BA {l}
 u_t-\Delta(|u|^{m-1}u)+h(t)|u|^{q-1}u=0\quad\mbox {in }Q^T
\\[2mm]
u(x,0)=u_{0,k}(x)=M_k^{\frac1{m+1}}k^{-\frac{mN}{m+1}}\delta_k(x)\qquad
 x\in\mathbb{R}^N,
\EA\right.\end{equation}
where $\delta_k$ is as in \rth{locth1}. Let us rewrite
problem \eqref{4.1} in the form
\begin{equation}\label{4.3}\left\{\BA {l}
 (|v|^{p-1}v)_t-\Delta v +h(t)|v|^{g-1}v=0,\quad\mbox {in }Q^T\\[2mm]
v=v_{k}=|u|^{m-1}u,\ p=1/m,\ g=q/m
\\[2mm]
|v(x,0)|^{p-1}v(x,0)=|v_{0,k}|^{p-1}v_{0,k}:=u_{0,k}(x)=M_k^{\frac
p{p+1}}k^{-\frac N{p+1}}\delta_k(x).
\EA\right.\end{equation}
Without loss of generality we may suppose
\begin{equation}\label{4.5}
\| \delta_k(x)
\|_{L_{\frac{p+1}p}(\mathbb{R}^N)}^{\frac{p+1}p}=\int_{\mathbb{R}^N}|\delta_k(x)|^{\frac{p+1}p}\,
dx\leq c_0k^{\frac Np}\qquad \forall\, k\in\mathbb{N}.
\end{equation}
Now sequence $\{M_k\}$ is such that
\begin{equation}\label{4.6}
M_k^{\frac p{p+1}}k^{-\frac N{p+1}}\to \infty\quad\text{as}\quad
k\to\infty.
\end{equation}
{\it Step 1. The local energy framework}.
Consider the following energy functions
\begin{equation}\label{4.7}
I_1(\tau)=\int_{Q_r}|\nabla_xv|^2\,dxdt,\qquad
I_2(\tau)=\int_{Q_r}h(t)|v|^{g+1}\,dxdt,\qquad
I_3(\tau)=\int_{Q_r}|v|^{p+1}\,dxdt.
\end{equation}
Analogously to \eqref{2.6} we deduce the inequality
\begin{multline}\label{4.8}
\int_{\mathbb{R}^N}|v(x,T)|^{p+1}\,dx+I_1(r)+I_2(r)+I_3(r)\leq
c\tau^{\frac{N(g-p)}{g+1}}h(r)^{-\frac{p+1}{g+1}}(-I'_2(r))^{\frac{p+1}{g+1}}
+c\int_{\Omega(\tau)}|v(x,r)|^{p+1}\,dx\\ \forall\, \tau>0, \
\forall\, r: 0<r<T.
\end{multline}
This inequality will control the spreading of energy with respect to the
 $r$-variable (the time direction). As to vanishing property of energy
in variable $\tau$, we will use the finite speed propagation of
support property for porous media equation with slow diffusion. In
the domain $Q^{(r)}(\tau)$ we will use the energy function
$E_1(r,\tau)=\int_{Q^{(r)}(\tau)}|\nabla_xv|^2\,dxdt$ from
\eqref{2.8}. Since $\text{supp}\,
v(\cdot,0)=\text{supp}\,v_k(\cdot,0)=\text{supp}\,v_{0,k}=\{x:|x|<k^{-1}\}$,
multiplying equation \eqref{4.3} on $v(x,t)$ and integrating in
the domain $Q^{(r)}(\tau),\ \tau\geq k^{-1}$, we obtain after
simple computations (see, for example \cite{DV,An}) the following
differential inequality
\begin{equation}\label{4.9}
\int_{\Omega(\tau)}|v(x,r)|^{p+1}\,dx+E_1(r,\tau)\leq
cr^{\frac{(p+1)(1-\theta_1)}{p+1-(1-\theta_1)(1-p)}}\Big(
-\frac{d}{d\tau}E_1(r,\tau)
\Big)^{\frac{p+1}{p+1-(1-\theta_1)(1-p)}},
\end{equation}
$$\forall\,\tau\geq
k^{-1},\ \forall\,r>0\quad \mbox {where }\;
\theta_1=\frac{N(1-p)+(p+1)}{N(1-p)+2(p+1)},\
1-\theta_1=\frac{p+1}{N(1-p)+2(p+1)}.
$$
Solving this inequality and keeping in mind that $E_1(r,\tau)\geq0\
\forall\,r>0,\ \forall\,\tau>0$, we deduce easily
\begin{equation}\label{4.10}
v(x,r)\equiv0\quad\forall\,
x:|x|>k^{-1}+c_0r^{1-\theta_1}E_1(r,k^{-1})^{\frac{(1-\theta_1)(1-p)}{1+p}}:=k^{-1}+c_0\chi(r),\quad\forall\,r>0.
\end{equation}
Here the constant $c_0>0$ depends on the parameters of the problem under
consideration, but do not on $r$ and $k$. Analogously to
\eqref{2.19} we deduce the following global {\it a priori} estimate
\begin{equation}\label{4.11}
\int_{Q^{(r)}}(|\nabla_xv|^2+r^{-1}|v|^{p+1}+h(t)|v|^{g+1})\,dxdt\leq
c\,\| v_{0,k} \|_{L_{p+1}(\mathbb{R}^N)}^{p+1}.
\end{equation}
Thus, due to \eqref{4.3}--\eqref{4.6}, it follows from
\eqref{4.11}
\begin{equation}\label{4.12}
E_1(r,0)\leq cM_k\qquad\forall\,r>0.
\end{equation}
Next we come back to the inequality \eqref{4.8}. Due to
\eqref{4.10} it ensues from \eqref{4.8} the inequality
\begin{equation}\label{4.12*}
I_1(r)+I_2(r)+I_3(r)\leq
c(k^{-1}+\chi(r))^{\frac{N(g-p)}{g+1}}h(r)^{-\frac{p+1}{g+1}}(-I'_2(r))^{\frac{p+1}{g+1}}\qquad
\forall\,r>0.
\end{equation}
Remark that due to \eqref{4.12} we have
\begin{equation}\label{4.13}
\chi(r)\leq
c_1r^{1-\theta_1}\,M_k^{\frac{(1-\theta_1)(1-p)}{1+p}}.
\end{equation}
{\it Step 2. The first round of computations}. Now we have to define $\tau_k,\ r_k $. First we impose the relation
\begin{equation}\label{4.14}
\tau_k\geq
c_1r_k^{1-\theta_1}\,M_k^{\frac{(1-\theta_1)(1-p)}{1+p}}, \qquad
c_1\ \text{is from (\ref{4.13})}.
\end{equation}
Then \eqref{4.12*} yields to
\begin{equation}\label{4.15}
I(r):=I_1(r)+I_2(r)+I_3(r)\leq
c(k^{-1}+\tau_k)^{\frac{N(g-p)}{g+1}}h(r)^{-\frac{p+1}{g+1}}(-I'(r))^{\frac{p+1}{g+1}}\qquad
\forall\,r:0<r<r_k.
\end{equation}
Solving this differential inequality we get the estimate
\begin{equation}\label{4.16}
I(r)\leq\frac{c(k^{-1}+\tau_k)^N}{\Big(\int_0^rh(s)\,ds\Big)^{\frac{p+1}{g-p}}}
\qquad\forall\,r:0<r<r_k.
\end{equation}
Remember that the function $h(s)$ has the form
$
h(s)=s^{(g-1)/(1-p)}\omega(s)^{-1},
$
therefore estimate \eqref{4.16} yields to
\begin{equation}\label{4.17}
I(r)\leq\frac{c_2\omega(r)^{\frac{p+1}{g-p}}(k^{-1}+\tau_k)^N}{r^{\frac{p+1}{1-p}}}
\qquad\forall\,r:0<r\leq r_k.
\end{equation}
Thus, as second relation, which defines our pair $\tau_k,\ r_k$,
we suppose the condition
\begin{equation}\label{4.18}
\frac{c_2\omega(r_k)^{\frac{p+1}{g-p}}(k^{-1}+\tau_k)^N}{r_k^{\frac{p+1}{1-p}}}\leq
cM_{k-1},\quad c \text{ is from (\ref{4.12})}.
\end{equation}
Moreover, we will find the pair $\tau_k,\ r_k$ such that
 the following property holds
\begin{equation}\label{4.19}
k^{-1}+\tau_k\leq1.
\end{equation}
Then the next inequality is a sufficient condition for validity
of \eqref{4.18}:
\begin{equation}\label{4.20}
c_2\omega(r_k)^{\frac{p+1}{g-p}}r_k^{-\frac{p+1}{1-p}}\leq
cM_{k-1},\quad c \text{ is from (\ref{4.12})},
\end{equation}
and we can define $r_k$ by equality
\begin{equation}\label{4.21}
r_k:=\Big(\frac{c_2}{c}\Big)^{\frac{1-p}{p+1}}\omega(r_k)^{\frac{1-p}{g-p}}M_{k-1}^{-\frac{1-p}{p+1}}.
\end{equation}
Now we have to choose the sequence $\{M_k\}$. Namely, we set
\begin{equation}\label{4.22}
M_k:=e^{k}\qquad\forall\,k\in\mathbb{N},
\end{equation}
and we define $\tau_k$, in accordance with assumption
\eqref{4.14}, by
\begin{equation}\label{4.23}
\tau_k=c_1r_k^{1-\theta_1}M_k^{\frac{(1-\theta_1)(1-p)}{1+p}},\quad
c_1\text{ is from (\ref{4.13})}.
\end{equation}
Further, due to \eqref{4.21} an \eqref{4.22}, it follows from
\eqref{4.23},
\begin{multline}\label{4.24}
\tau_k=c_1(r_k^{p+1}M_k^{1-p})^{\frac1{N(1-p)+2(p+1)}}= c_1\Big[
\Big( \frac{c_2}c
\Big)^{1-p}\omega(r_k)^{\frac{(1-p)(p+1)}{g-p}}M_{k-1}^{-(1-p)}M_k^{1-p}
\Big]^{\frac1{N(1-p)+2(p+1)}}\\=c_1\Big( \frac{ec_2}c
\Big)^{\frac{(1-\theta_1)(1-p)}{1+p}}\omega(r_k)^{S},
\end{multline}
where
$S=\frac{(1-\theta_1)(1-p)}{g-p}=\frac{(1-p)(p+1)}{(g-p)[N(1-p)+2(p+1)]}$.
From definition \eqref{4.21} and because of \eqref{4.22} and \eqref{2.35},
there holds
\begin{equation}\label{4.25}
r_k\leq\Big( \frac{c_2}c
\Big)^{\frac{1-p}{p+1}}\omega_0^{\frac{1-p}{g-p}}\exp\Big(
-\frac{1-p}{p+1}(k-1) \Big):=c_3\exp\Big( -\frac{1-p}{p+1}\,k
\Big),
\end{equation}
and $r_k\to0$ as $k\to\infty$. Therefore, since $\omega(s)\to0$ as
$s\to0$, it follows from \eqref{4.24} that $\tau_k\to0$ as
$k\to\infty$. Consequently we can suppose $k$ so large that condition
\eqref{4.19} is satisfied. Thus, we have pair $(\tau_k,\
r_k)$ for large $k\in\mathbb{N}$. \\

 \noindent{\it Step 3. The second round of computations}.  As a starting global {\it a priori} estimate of solution we
will use now, instead of \eqref{4.11}, \eqref{4.12}, the following
estimate
\begin{equation}\label{4.26}
I_1(r_k)=\int_{\substack{\{t\geq r_k,\\ x\in\mathbb{R}^N\}}}
|\nabla_xv|^2\,dxdt\leq I(r_k)\leq cM_{k-1},
\end{equation}
which follows from \eqref{4.17}, due to definition \eqref{4.18},
\eqref{4.21} of $r_k$.  Using property \eqref{4.10}, estimate
\eqref{4.13} and property \eqref{4.14}, it ensues from
\eqref{4.26}
\begin{equation}\label{4.27}
E_1(r,k^{-1}+\tau_k)\leq I_1(r)\leq I_1(r_k)<cM_{k-1}\qquad
\forall\,r\geq r_k.
\end{equation}
Since $v(x,r_k)=0\ \forall\,x:|x|\geq k^{-1}+\tau_k$ we deduce
similarly to \eqref{4.9}
\begin{multline}\label{4.28}
\int_{\Omega(\tau)}|v(x,r_k+r)|^{p+1}\,dx+E_1(r_k+r,k^{-1}+\tau_k+\tau)\leq
cr^{\frac{(p+1)(1-\theta_1)}{(p+1)-(1-\theta_1)(1-p)}}\\\times\Big(-\frac
d{d\tau}E_1(r_k+r,k^{-1}+\tau_k+\tau)
\Big)^{\frac{p+1}{p+1-(1-\theta_1)(1-p)}}\qquad \forall\,r>0,\
\forall\,\tau>0.
\end{multline}
Solving this differential inequality, we obtain
\begin{equation}\label{4.29}
v(x,r_k+r)\equiv0\quad \forall\, x:|x|\geq
k^{-1}+\tau_k+c_0\chi_1(r),
\end{equation}
where
$\chi_1(r):=r^{1-\theta_1}E_1(r_k+r,k^{-1}+\tau_k)^{\frac{(1-\theta_1)(1-p)}{1+p}}\
\forall\,r\geq 0. $ But \eqref{4.27} implies
\begin{equation}\label{4.30}
\chi_1(r)\leq
c_1r^{1-\theta_1}M_{k-1}^{\frac{(1-\theta_1)(1-p)}{1+p}}.
\end{equation}
Now we define $\tau_{k-1},\ r_{k-1}$.
In the same way as \eqref{4.14} we impose
\begin{equation}\label{4.31}
\tau_{k-1}\geq
c_1r_{k-1}^{1-\theta_1}M_{k-1}^{\frac{(1-\theta_1)(1-p)}{1+p}}.
\end{equation}
Similarly to \eqref{4.15}--\eqref{4.17} we deduce
\begin{equation}\label{4.32}
I(r)\leq\frac{c_2\omega(r)^{\frac{p+1}{g-p}}(k^{-1}+\tau_k+\tau_{k-1})^N}{r^{\frac{p+1}{1-p}}}\qquad
\forall\,r:0<r\leq r_k+r_{k-1}.
\end{equation}
The second relation for defining the pair $\tau_{k-1},\
r_{k-1}$ is analogous to\eqref{4.18}
\begin{equation}\label{4.33}
\frac{c_2\omega(r_k+r_{k-1})^{\frac{p+1}{g-p}}(k^{-1}
+\tau_k+\tau_{k-1})^N}{(r_k+r_{k-1})^{\frac{p+1}{1-p}}}\leq
cM_{k-2},\quad c\text{ is from (\ref{4.12})}.
\end{equation}
Supposing that
\begin{equation}\label{4.34}
k^{-1} +\tau_k+\tau_{k-1}\leq1,
\end{equation}
we can define $r_{k-1}$ by the following analogue of \eqref{4.21}
\begin{equation}\label{4.35}
r_k+r_{k-1}:=\Big( \frac{c_2}c
\Big)^{\frac{1-p}{p+1}}\omega(r_k+r_{k-1})^{\frac{1-p}{g-p}}M_{k-2}^{-\frac{1-p}{p+1}}.
\end{equation}
And in accordance with \eqref{4.31} let us define $\tau_{k-1}$ by
\begin{equation}\label{4.36}
\tau_{k-1}=
c_1r_{k-1}^{1-\theta_1}M_{k-1}^{\frac{(1-\theta_1)(1-p)}{1+p}}.
\end{equation}
Due to \eqref{4.35} we have
\begin{equation*}\label{4.37}\BA {l}
\tau_{k-1}\leq c_1\big[ (r_k+r_{k-1})^{p+1}M_{k-1}^{1-p}
\big]^{\frac1{N(1-p)+2(p+1)}}\\[3mm]
\phantom{----------}\leq c_1\Big[
\Big(\myfrac{c_2}{c}\Big)^{1-p}\omega(r_k+r_{k-1})^{\frac{(1-p)(p+1)}{g-p}}M_{k-2}^{-(1-p)}M_{k-1}^{1-p}
\Big]^{\frac1{N(1-p)+2(p+1)}}\\[3mm]
\phantom{----------}=c_1\Big( \myfrac{ec_2}c
\Big)^{\frac{(1-\theta_1)(1-p)}{1+p}}\omega(r_k+r_{k-1})^{S},
\EA\end{equation*} where $S$ is from (\ref{4.24}). Notice that,
due to \eqref{4.32}, \eqref{4.33}, we have also
\begin{equation}\label{4.38}
I_1(r_k+r_{k-1})\leq I(r_k+r_{k-1})\leq cM_{k-2},
\end{equation}
and, analogously to \eqref{4.27},
\begin{equation}\label{4.39}
E_1(r,k^{-1}+\tau_k+\tau_{k-1})\leq I_1(r)\leq
I_1(r_k+r_{k-1})\leq cM_{k-2}\qquad \forall\, r\geq r_k+r_{k-1}.
\end{equation}
 \noindent{\it Step 4. Completion of the proof}.
Estimates \eqref{4.38}, \eqref{4.39} we can use instead of
\eqref{4.26}, \eqref{4.27} for third round of computations. After
$j$ such rounds we deduce that
\begin{align}\label{4.40}
&I_1\bigg( \sum_{i=0}^j r_{k-i} \bigg)\leq I\bigg( \sum_{i=0}^j
r_{k-i} \bigg)\leq c M_{k-j},\\
\label{4.41} & E_1\bigg(  r,k^{-1}+\sum_{i=0}^j\tau_{k-i}
\bigg)\leq I_1(r)\leq I_1\bigg(  \sum_{i=0}^j r_{k-i} \bigg)\leq
cM_{k-j}\qquad \forall\,r\geq\sum_{i=0}^j r_{k-i},
\end{align}
where
\begin{equation}\label{4.42}
\tau_{k-i}\leq
c_1\Big(\frac{ec_2}c\Big)^{\frac{(1-\theta_1)(1-p)}{1+p}}
\omega\bigg(\sum_{l=0}^ir_{k-l}\bigg)^{S},
\end{equation}
with the same $S$ as in \eqref{4.24}, and
\begin{equation}\label{4.43}
\sum_{l=0}^ir_{k-l}=\Big(\frac{c_2}c\Big)^{\frac{1-p}{p+1}}
\omega\bigg(\sum_{l=0}^ir_{k-l}\bigg)^{\frac{1-p}{g-p}}M_{k-i-1}^{-\frac{1-p}{p+1}}.
\end{equation}
Estimates \eqref{4.40} will remain true as long as the
following analogue of relation \eqref{4.34} is valid
\begin{equation*}\label{4.44}
k^{-1}+\sum_{i=0}^j\tau_{k-i}\leq1.
\end{equation*}
Now we will check this condition. Due to \eqref{2.25}, it follows
from \eqref{4.43}
\begin{equation*}\label{4.45}
\sum_{l=0}^ir_{k-l}\leq\Big(\frac{c_2}c\Big)^{\frac{1-p}{p+1}}\omega_0^{\frac{1-p}{g-p}}
M_{k-i-1}^{-\frac{1-p}{p+1}}:=CM_{k-i-1}^{-\frac{1-p}{p+1}}=C\exp\Big(-\frac{1-p}{p+1}(k-i-1)\Big).
\end{equation*}
Therefore, from \eqref{4.42}, it follows
\begin{equation*}\label{4.46}\BA {l}
\tau_{k-i}\leq c_1\Big( \myfrac{ec_2}c
\Big)^{\frac{(1-\theta_1)(1-p)}{1+p}}\omega\Big(
C\exp\Big(-\frac{(1-p)(k-i-1)}{p+1}\Big)
\Big)^{S}:=C_1\Big[\omega\Big(
C\exp\Big(-\frac{(1-p)(k-i-1)}{p+1}\Big) \Big)\Big]^S.
\EA\end{equation*} Thus we have, using in particular the
monotonicity of function $\omega(s)$,
\begin{multline}\label{4.47}
\displaystyle \sum_{i=0}^{j}\tau_{k-i}\leq C_1\sum_{i=0}^{j} \Big[
\omega\Big(C\exp\Big(-\frac{(1-p)(k-i-1)}{p+1}\Big) \Big)
\Big]^S\\
\leq C_1\myint{k-j-1}{k}\Big[ \omega\Big(
C\exp\Big(-\frac{(1-p)s}{p+1}\Big) \Big) \Big]^S\,ds
 =\myfrac{C_1(p+1)}{1-p}
\myint{A_1}{A_2}
\myfrac{\omega(s)^S}s\,ds,\\
A_1=C\exp\big(-\frac{1-p}{p+1}k\big),\quad
A_2=C\exp\big[-\frac{1-p}{p+1}(k-j-1)\big].
\end{multline}
 Due to condition
\eqref{intr6*} and estimate \eqref{4.47} we can find
$k_0\in\mathbb{N}$, which depends on parameters of problem under
consideration, but does not depend on $k\in\mathbb{N}$, such that
$$
\sum_{i=0}^{k-k_0}\tau_{k-i}+k^{-1}\leq1\qquad\forall\,k\in\mathbb{N}.
$$
At end, our estimates \eqref{4.40}--\eqref{4.43} are true for all
$j\leq k-k_0$. Therefore the proof of \rth{pmth2}  follows
from estimates \eqref{4.40}--\eqref{4.43}, in the same way as
\rth{locth1} from estimates
\eqref{2.70}--\eqref{2.72}.\qeda

\section{The fast diffusion equation with absorption}
\setcounter{equation}{0}

When $(1-2/N)_{+}<m<1$, it is known that the mere fast diffusion equation
\begin {equation}\label {FDE}
\prt_{t}v-\Gd v^m=0\quad\mbox {in }\;Q^\infty
\end{equation}
admits a particular fundamental positive solution with initial data $k\gd_{0}$ ($k>0$) called the Barenblatt -Zeld'dovich-Kompaneets
solution, expressed by
\begin{equation}\label{barb3}
B_k(x,t)=t^{-\ell}\left(C_k+\myfrac{(1-m)\ell}{2mN}\myfrac{\abs x^2}{t^{2\ell/N}}\right)
^{-1/(1-m)},
\end {equation}
where $\ell$ and $C_{k}$ are given in (\ref{barb2}). The main feature of this expression is that
$\lim_{k\to\infty}C_{k}=0$, therefore
\begin{equation}\label{barb4}
\lim_{k\to\infty}B_k(x,t)=W(x,t):=C_{*}\left(\myfrac {t}{\abs x^2}\right)
^{1/1-m)},
\end {equation}
where
$$C_{*}=\left(\myfrac{(1-m)^3}{2m(mN+2-N)}\right)^{1/(1-m)}.
$$
This solution has a persisting singularity and is called a razor blade \cite {VV}. It has also the property that
$$\lim_{t\to 0}W(x,t)=0\forevery x\neq 0.
$$
This phenomenon is at the origin of the work of Chasseigne and V\`azquez on extended solutions of the fast diffusion equation \cite {CV}.
Concerning problem (\ref{pme-k}), \rprop {fundpm} is still valid provided $m>(1+2/N)_{+}$. We shall denote by $u=u_{k}$ the solutions of (\ref{pme-k}). Furthermore estimate (\ref{max1}) holds. Combining this with the fact that the $B_{k}$ are super solutions for the $u_{k}$, we derive the following
\bth{fdeth} Assume $(1-2/N)_{+}<m<1$ and $h\in C(0,\infty)$ is positive. Assume also that
 \eqref{Bar-q1} holds. Then $u_{\infty}:=\lim_{k\to\infty}u_{k}$ has a point-wise singularity at $(0,0)$ and the following estimate is verified
\begin{equation}\label{barb5}
u_{\infty}(x,t)\leq\min\left\{C_{*}t^{-\ell}
\left(\myfrac{\abs x^2}{t^{2\ell/N}}\right)^{-1/(1-m)},\left((q-1)\myint {0}{t}h(s)\,ds\right)^{-1/(q-1)}\right\}
\end {equation}
\es

\noindent \Remark The profile of $u_{\infty}$ near $(x,t)=(0,0)$ is completely unknown. In particular a very chalenging question could be to give precise estimates on the quantity $\min\left\{W(x,t),U_{h}(t)\right\}-u_{\infty}(x,t)$.

\enddocument